\documentclass[10pt,a4paper]{article}

\usepackage{amsmath,amsthm,amssymb,graphics,setspace}
\usepackage{calc}
\usepackage[dvips]{color} 
\usepackage{bm}
\newtheorem{thm}{\hspace*{20pt}Theorem}
\newtheorem{lem}{\hspace*{20pt}Lemma}[section]
\newtheorem{pro}{\hspace*{20pt}Proposition}

\newcounter{constant}[section]
\newcommand{\const}{\ifnum \theconstant > 0 \stepcounter{constant}\theconstant
                                             \else \setcounter{constant}{1}\theconstant \fi }
\newcounter{subconstant} 
\newcounter{submioneconstant}

\newcounter{subconstanta}

\newcounter{submitwoconstant}

\newcommand{\suba}{\thesubconstanta} 
\newcommand{\insa}{\setcounter{subconstanta}{\value{constant}}}  
\newcounter{subconstantb} 
\newcommand{\subb}{\thesubconstantb} 
\newcommand{\insb}{\setcounter{subconstantb}{\value{constant}}}  
\newcounter{subconstantc} 
\newcommand{\subc}{\thesubconstantc} 
\newcommand{\insc}{\setcounter{subconstantc}{\value{constant}}}  
\newcounter{subconstantd} 
\newcommand{\subd}{\thesubconstantd} 
\newcommand{\insd}{\setcounter{subconstantd}{\value{constant}}}  
\newcounter{subconstante} 
\newcommand{\sube}{\thesubconstante} 
\newcommand{\inse}{\setcounter{subconstante}{\value{constant}}}  
\newcounter{subconstantf} 
\newcommand{\subf}{\thesubconstantf} 
\newcommand{\insf}{\setcounter{subconstantf}{\value{constant}}}  
\newcounter{subconstantg} 
\newcommand{\subg}{\thesubconstantg} 
\newcommand{\insg}{\setcounter{subconstantg}{\value{constant}}}  
\newcounter{subconstanth} 
\newcommand{\subh}{\thesubconstanth} 
\newcommand{\insh}{\setcounter{subconstanth}{\value{constant}}}  
\newcounter{subconstanti}

\newcounter{subconstantj}

\newcounter{subconstantk}

\newcounter{subconstantl}

\newcounter{subconstantm}

\newcounter{subconstantn}

\newcounter{subconstanto}

\newcounter{subconstantp} 
\newcommand{\subp}{\thesubconstantp} 
\newcommand{\insp}{\setcounter{subconstantp}{\value{constant}}}  
\newcounter{subconstantq} 
\newcommand{\subq}{\thesubconstantq} 
\newcommand{\insq}{\setcounter{subconstantq}{\value{constant}}}  
\newcounter{subconstantr} 
\newcommand{\subr}{\thesubconstantr} 
\newcommand{\insr}{\setcounter{subconstantr}{\value{constant}}}  
\newcounter{subconstants}

\newcounter{subconstantt}

\newcounter{subconstantu}

\newcounter{subconstantv}

\newcounter{subconstantw}

\newcounter{subconstantx} 
\newcommand{\subx}{\thesubconstantx} 
\newcommand{\insx}{\setcounter{subconstantx}{\value{constant}}} 
\newcounter{subconstanty} 
\newcommand{\suby}{\thesubconstanty} 
\newcommand{\insy}{\setcounter{subconstanty}{\value{constant}}}   
\newcounter{subconstantz} 
\newcommand{\subz}{\thesubconstantz} 
\newcommand{\insz}{\setcounter{subconstantz}{\value{constant}}}

\newcommand{\Rd}{{\bf R}^2}
\newcommand{\lk}{\lambda_k}
\newcommand{\vk}{v_k}

\newcommand{\tvk}{\tilde{v}_k}
\newcommand{\wk}{w_k}
\newcommand{\wkn}{w_{k,1}}
\newcommand{\twkn}{\tilde{w}_{k,1}}
\newcommand{\wknn}{w_{k,2}}
\newcommand{\psik}{\psi_k}
\newcommand{\tpsik}{\tilde{\psi}_k}
\newcommand{\hk}{h_k}
\newcommand{\thk}{\tilde{h}_k}

\newcommand{\zkd}{z_{k,1}^\delta}  
\newcommand{\zkkd}{z_{k,2}^\delta}  
\newcommand{\lzero}{\lambda_0}
\newcommand{\rk}{r_k}
\newcommand{\Gg}{{\cal G}}
\newcommand{\cS}{{\cal S}}
\newcommand{\cSp}{{\cal S}_{+}}
\newcommand{\cSm}{{\cal S}_{-}}
\newcommand{\cSpm}{{\cal S}_{\pm}}
\newcommand{\cSmp}{{\cal S}_{\mp}}
\newcommand{\xk}{x_k}
\newcommand{\supp}{{\rm supp}\ }
\newcommand{\diam}{{\rm diam}}

\makeatletter
\@addtoreset{equation}{section}

\makeatother

\begin{document}

\title{Blowup behavior for a degenerate elliptic 
$\sinh$-Poisson equation with variable intensities}

\author{Tonia Ricciardi\footnotemark[1] \ \ \ Ryo Takahashi \footnotemark[2]}

\footnotetext[1]{
Dipartimento di Matematica e Applicazioni ``R. Caccioppoli''
Universit\`{a} di Napoli Federico II, Via Cintia, 80126 Napoli, Italy
(E-mail: {\it tonia.ricciardi@unina.it})
}
\footnotetext[2]{
Division of Mathematical Science, Department of Systems Innovation, 
Graduate School of Engineering Science, Osaka University, 
Machikaneyamacho 1-3, Toyonakashi, 560-8531, Japan. 
(E-mail: {\it r-takaha@sigmath.es.osaka-u.ac.jp})
}
\date{\today}
\maketitle
\begin{abstract}
In this paper, we provide a complete blow-up picture for solution sequences 
to an elliptic sinh-Poisson equation with variable intensities
arising in the context of the statistical mechanics description of two-dimensional turbulence,
as initiated by Onsager. 
The vortex intensities are described in terms of a probability measure $\mathcal P$
defined on the interval $[-1,1]$.
Under Dirichlet boundary conditions we establish 
the exclusion of boundary blowup points,
we show that the concentration mass does not have
residual $L^1$-terms (``residual vanishing") and 
we determine the location of blowup points
in terms of Kirchhoff's Hamiltonian.
We allow $\mathcal P$
to be a general Borel measure,
which could be ``degenerate" in the sense that 
$\mathcal P(\{\alpha_-^*\})=0=\mathcal P(\{\alpha_+^*\})$, 
where $\alpha_-^*=\min\mathrm{supp}\mathcal P$ 
and $\alpha_+^*=\max\mathrm{supp}\mathcal P$.
Our main results are new for the standard sinh-Poisson equation as well.
\end{abstract}

\section{Introduction}\label{sec:intro}

Since Onsager's pioneering work \cite{onsager} in 1949, 
the statistical mechanics description of stable, large scale vortices has attracted
the attention of many physicists and mathematicians, and is still of central interest in fluid mechanics 
\cite{BouchetVenaille, Chavanis}.
In particular, several mean field equations have been proposed to describe 
two-dimensional stationary Euler flows 
with a large number of point vortices. 
\par
In this paper we are concerned with the following mean field equation derived by C.~Neri in \cite{neri}
under the ``stochastic" assumption that the vortex intensities and orientations are 
independent identically distributed random
variables with probability distribution $\mathcal P(d\alpha)$, $\alpha\in[-1,1]$:
\begin{equation}\label{eqn:neri}
\left\{
\begin{aligned}
\displaystyle 
-\Delta v=&\lambda \int_{[-1,1]} \frac{\alpha e^{\alpha v}}{\iint_{[-1,1]\times \Omega} 
e^{\alpha'v}\,{\cal P}(d\alpha')dx}\,{\cal P}(d\alpha) && \mbox{in $\Omega$} \\
\displaystyle 
v=&0 && \mbox{on $\partial\Omega$}. 
\end{aligned}
\right.
\end{equation} 
Here, $v$ denotes the stream function of a turbulent Euler flow, 
$\Omega\subset{\bf R}^2$ is a bounded domain with smooth boundary $\partial\Omega$, 
$\lambda>0$ is a constant related to the inverse temperature. 
We further assume that
${\cal P}\in{\cal M}([-1,1])$ is a Borel probability measure defined on the interval $[-1,1]$, 
where ${\cal M}([-1,1])$ denotes the space of measures on
$[-1,1]$. 
\par
If $\mathcal P=\delta_{+1}$, corresponding to the case where all vortices have the same intensity and orientation,
equation~\eqref{eqn:neri} reduces to the Liouville type equation
\begin{equation}\label{eqn:meanfield}
\left\{
\begin{aligned}
-\Delta v=&\lambda\frac{e^v}{\int_\Omega e^v\,dx} && \mbox{in $\Omega$} \\
v=&0 && \mbox{on $\partial\Omega$}
\end{aligned}
\right.
\end{equation}
whose properties are by now well understood,
see, e.g., \cite{CSLin, suz08} and the references therein. 
If ${\cal P}=(\delta_{+1}+\delta_{-1})/2$,
corresponding to the case
where the point vortices have the same intensity and variable orientation,
equation~\eqref{eqn:neri} reduces to the sinh-Poisson type problem:
\begin{equation}\label{eqn:sinhGordon}
\left\{
\begin{aligned}
-\Delta v=&\frac{\lambda}{2}\frac{e^v-e^{-v}}{\int_\Omega(e^v+e^{-v})\,dx} && \mbox{in $\Omega$} \\
v=&0 && \mbox{on $\partial\Omega$}
\end{aligned}
\right.
\end{equation}
Equation~\eqref{eqn:sinhGordon} is also related to the study of constant mean curvature surfaces and has
received a considerable attention, see,
e.g., \cite{BartolucciPistoia, grossipistoia, jostwangyezhou, spruck}
and the references therein.
Our results for \eqref{eqn:neri} will yield new results for \eqref{eqn:sinhGordon}
as well. 
\par
It is useful to mention that another mean field equation with probability measures formally similar to \eqref{eqn:neri}
was derived under a ``deterministic" assumption
on the vortex intensities in \cite{ss2008}, see also Onsager's handwritten note in \cite{es2006}: 
\begin{equation}\label{eqn:ss}
\left\{
\begin{array}{ll}
\displaystyle 
-\Delta v=\lambda \int_{[-1,1]} \frac{\alpha e^{\alpha v}}{\int_{\Omega} e^{\alpha v} dx} {\cal P}(d\alpha) & \mbox{in $\Omega$} \\
\displaystyle 
v=0 & \mbox{on $\partial\Omega$}. 
\end{array}
\right.
\end{equation} 
Equation~\eqref{eqn:ss} also reduces to \eqref{eqn:meanfield} when ${\cal P}=\delta_{+1}$.
However, if ${\cal P}=(\delta_{+1}+\delta_{-1})/2$,
equation~\eqref{eqn:ss} reduces to
\begin{equation}\label{eqn:pl}
\left\{
\begin{aligned}
-\Delta v=&\frac{\lambda}{2}\left(\frac{e^v}{\int_\Omega e^v\,dx}-\frac{e^{-v}}{\int_\Omega e^{-v}\,dx}\right)&& \mbox{in $\Omega$} \\
v=&0 && \mbox{on $\partial\Omega$}
\end{aligned}
\right.
\end{equation}
which is evidently different from \eqref{eqn:sinhGordon}.
Results for equation \eqref{eqn:pl} may be found in \cite{EspositoWei, os2006}
and the references therein.
\par
It is therefore a natural question to seek common properties between \eqref{eqn:neri} and \eqref{eqn:ss}
as well as different behaviors,
which could in principle provide a criterion to select a ``more suitable" model.
Several results in this direction were obtained in \cite{ors, ricciardisuzuki, rz2012, rz2014, rz2014ANS}.
In particular, with the aim of studying common properties of solution sequences on compact surfaces, a very general equation
containing \eqref{eqn:neri} and \eqref{eqn:ss} as special cases
was considered in \cite{rz2012}. 
Roughly speaking, it was shown that the basic Brezis-Merle type blow-up
alternatives holds true for both models, see Proposition \ref{prop:rz2012} below.
On the other hand, 
consideration of the optimal Moser-Trudinger type inequalities
associated to \eqref{eqn:neri} and \eqref{eqn:ss} 
emphasized significantly different properties
between \eqref{eqn:neri} and \eqref{eqn:ss},
see \cite{ricciardisuzuki, rz2012, stz}.
\par
\par
Our aim in this article is to complete the blow-up analysis for solution
sequences to \eqref{eqn:neri} initiated in \cite{rz2012, rz2014}.
We shall first of all show that under Dirichlet boundary conditions
blow-up cannot occur on the boundary.
We notice that the exclusion of boundary blowup points (see Theorem \ref{thm:main}-{\bf II-(i)}
for the precise statement)
is not straightforward for general cases of $\mathcal P$, 
although it is readily derived from an estimate in \cite{gnn} for the one-sided case, 
that is, the case where $\supp{\cal P}\subset [0,1]$ or $\supp{\cal P}\subset [-1,0]$. 
In a forthcoming paper \cite{rtzz-pre} we shall show that the exclusion
of boundary blow-up points holds true for a more general class of equations
including \eqref{eqn:neri} and \eqref{eqn:ss} with Dirichlet boundary conditions.
Then, we establish
the vanishing of the $L^1$-terms in the concentration mass limit
(``residual vanishing"), see Theorem \ref{thm:main}-{\bf II-(ii)}
or Section \ref{sec:rv}. This property was derived in \cite{rz2014}
for $\mathcal P$ satisfying 
$\mathrm{supp}\mathcal P\subset[0,1]$ and
$\mathcal P(\{1\})>0$, but the case of a general Borel measure $\mathcal P$
was left open. 
We note that the residual vanishing property is specific
to \eqref{eqn:neri}, in the sense that it is known to be false for 
\eqref{eqn:ss} for some
special choices of $\mathcal P$, see \cite{stz}.
As a consequence of the residual vanishing property we shall 
locate the blow-up points
in terms of Kirchhoff's Hamiltonian \cite{BartschPistoia, clmp92, suz08}
\begin{equation}
\label{def:KRHamiltonian}
\mathcal H_N(x_1,\ldots,x_N)=\sum_{i=1}^Nr_i^2H(x_i,x_i)+\sum_{\stackrel{i,j=1}{i\neq j}}^Nr_ir_jG(x_i,x_j),
\end{equation}
where $r_i\in[-1,1]$, $i=1,\ldots,N$ denotes the vortex intensity of $x_i$
and where the sign of $r_i$ determines the vortex orientation. 
Here,
$G=G(x,y)$ and $H=H(x,y)$ are the Green's function and its regular part, respectively, that is, 
\begin{equation}\label{eqn:green-dom}
-\Delta G(x,y)=\delta_y\ \mbox{in $\Omega$},\quad G(x,y)=0\ \mbox{on $\partial\Omega$} 
\end{equation}
and
\begin{equation}\label{eqn:robin-dom}
H(x,y)=G(x,y)-\frac{1}{2\pi}\log\frac{1}{|x-y|}. 
\end{equation}
We recall that the starting point for the statistical mechanics derivation of mean field equations
for stationary flows with many vortices is given by Kirchhoff's point vortex model 
whose dynamics is governed by the Hamiltonian $\mathcal H_N$.
The mean field equations \eqref{eqn:neri}--\eqref{eqn:ss} are then derived by statistical mechanics 
arguments letting $N\to\infty$, along some ideas in \cite{clmp92}.
It is therefore expected that if the residual vanishing property holds, then solutions 
to \eqref{eqn:neri}--\eqref{eqn:ss} should concentrate at critical points for $\mathcal H_N$.
Theorem \ref{thm:lbp} will rigorously establish this fact for equation~\eqref{eqn:neri}.
\par
In order to state our  main results more precisely, we introduce some notation.
Let $(\lk,\vk)$ be a solution sequence to \eqref{eqn:neri}. 
We define the blowup sets:
\begin{align*}
&\cS_{\pm}=\{x_0\in\overline{\Omega}\ | \ \mbox{there exists $\xk\in\Omega$ 
such that  $\xk\rightarrow x_0$ and $\vk(\xk)\rightarrow\pm\infty$} \},\\ 
&\cS=\cSp\cup\cSm. 
\end{align*}
We point out that our definition of $\mathcal S$ allows the case $\mathcal S\cap\partial\Omega\neq\emptyset$.
We further define
\begin{equation}
\label{def:suppP}
\begin{aligned}
&\alpha_-^*=\min\mathrm{supp}\mathcal P,
&&\alpha_+^*=\max\mathrm{supp}\mathcal P
\end{aligned}
\end{equation}
where $\supp{\cal P}=\{\alpha\in[-1,1]\ | \ {\cal P}(N)>0\ \mbox{for any neighborhood $N$ of $\alpha$}\}$
denotes the support of $\mathcal P$.
For every $x_0\in\mathcal S$ we set
\begin{equation}\label{eqn:mq4}
\beta_\pm(x_0)=\begin{cases}
|\alpha_\pm^\ast|^{-1} & \mbox{if $x_0\in\cS_\pm$} \\
0 & \mbox{if $x_0\not\in\cS_\pm$}
\end{cases}.
\end{equation}
With this notation, we have the following.
\begin{thm}\label{thm:main}
Let $(\lk,\vk)$ be a solution sequence for \eqref{eqn:neri}.
Then, passing to a subsequence, we have the following alternatives. 
\begin{enumerate}
\item[{\bf (I)}] Compactness: $\displaystyle \limsup_{k\rightarrow\infty}\|\vk\|_\infty<+\infty$, that is, $\cS=\emptyset$. \\
Then, there exists $v\in H_0^1(\Omega)$ such that $\vk\rightarrow v$ in $H_0^1(\Omega)$ and $v$ is a solution of \eqref{eqn:neri}. 
\item[{\bf (II)}] Concentration: $\displaystyle \limsup_{k\rightarrow\infty}\|\vk\|_\infty=\infty$, that is, $\cS\neq\emptyset$. \\
Then, the following properties hold:
\begin{enumerate}
\item[{\bf (i)}][Exclusion of boundary blowup points]:
\begin{equation}
\label{eqn:exclusion} 
\mathcal S\cap\partial\Omega=\emptyset. 
\end{equation}
\item[{\bf (ii)}] It holds that $v_k \rightarrow v_0$ in $C_{loc}^2(\Omega\setminus\cS)$, 
where $G=G(x,y)$ is the Green function defined by \eqref{eqn:green-dom} and 
\begin{equation}\label{eqn:residual-vanishing}
v_0(x)=\sum_{x_0'\in\cS}(m_+(x_0')-m_-(x_0'))G(x,x_0')
\end{equation}
with $m_\pm(x_0')\geq 4\pi$ for every $x_0'\in\cS$.
\item[{\bf (iii)}][Mass relation]
\begin{equation}
(m_+(x_0)-m_-(x_0))^2=8\pi(\beta_+(x_0) m_+(x_0)+\beta_-(x_0)m_-(x_0))
 \label{eqn:mr} 
\end{equation}
for every $x_0\in\cS$, where $\beta_\pm(x_0)$ is defined in \eqref{eqn:mq4}.
\end{enumerate}
\end{enumerate}
\end{thm}
The mass relation \eqref{eqn:mr} was first noticed for \eqref{eqn:ss} in \cite{os2006}
in the special case where ${\cal P}$ is given by ${\cal P}=\tau\delta_{+1}+(1-\tau)\delta_{-1}$ ($\tau\in(0,1)$). 
It was then derived for \eqref{eqn:ss} with a general probability $\mathcal P$ in \cite{ors}.
In \cite{rz2012} a mass relation was established for a general equation including \eqref{eqn:neri}
and \eqref{eqn:ss} as special cases, see Proposition \ref{prop:rz2012} below.
\par
Our second results is concerned with the location of the blow-up points. 
\begin{thm}\label{thm:lbp}
Let $(\lk,\vk)$ be a solution sequence for \eqref{eqn:neri} and suppose the alternative 
{\bf (II)} in Theorem \ref{thm:main} occurs. 
Then, for every $x_0\in{\cal S}$, we have 
\begin{equation}\label{eqn:lbp-dom}
\nabla\left.\left[
H(x,x_0)+\sum_{x_0'\in{\cal S}\setminus\{x_0\}}\frac{(m_+(x_0')-m_-(x_0'))}{(m_+(x_0)-m_-(x_0))}G(x,x_0')
\right]\right|_{x=x_0}=0. 
\end{equation}
\end{thm}
We note that at a blow-up point we necessarily have $m_+(x_0)-m_-(x_0)\neq0$
in view of the mass relation Theorem \ref{thm:main}-{\bf II-(iii)},
so that the function in \eqref{eqn:lbp-dom} is well-defined.
Setting $\mathcal S=\{x_1,\ldots,x_m\}$, 
\eqref{eqn:lbp-dom} is equivalent to stating that $(x_1,\ldots,x_m)$
is a critical point for the Hamiltonian $\mathcal H_N$ defined in \eqref{def:KRHamiltonian}
with $N=m=\mathrm{card}\,\mathcal S$ and $r_i=m_+(x_i)-m_-(x_i)$,
$i=1,\ldots,m$.
\par
We observe that the properties stated in Theorem \ref{thm:main} and Theorem \ref{thm:lbp}
are essentially properties of the following ``model case"
\begin{equation}\label{eqn:model}
\left\{
\begin{aligned}
\displaystyle 
-\Delta v=&\frac{\lambda}{2}\, \frac{\alpha_+^* e^{\alpha_+^* v}+\alpha_-^*e^{\alpha_-^* v}}
{\int_\Omega(e^{\alpha_+^* v}+e^{\alpha_-^* v})}&& \mbox{in $\Omega$} \\
\displaystyle 
v=&0 && \mbox{on $\partial\Omega$}
\end{aligned}
\right.
\end{equation} 
corresponding to $\mathcal P=(\delta_{\alpha_+^*}+\delta_{\alpha_-^*})/2$
(recall that $\alpha_-^*\le0$).
This fact is evident in the  
``nondegenerate case" $\mathcal P(\{\alpha_+^*\})>0$ and  $\mathcal P(\{\alpha_-^*\})>0$;
technical care is needed in order to show that it holds true
in the  ``degenerate case" $\mathcal P(\{\alpha_+^*\})=\mathcal P(\{\alpha_-^*\})=0$.
We note that sign-changing nodal solutions for \eqref{eqn:model} blowing up at two distinct points
of $\Omega$ were recently constructed in \cite{pistoiaricciardi}.
\par
This paper is organized as follows.
In Section \ref{sec:pre}, we provide several preliminary lemmas.
In particular, we construct a convenient conformal mapping $X_0$ which will be useful
to ``straighten the boundary" $\partial\Omega$ near a point $x_0\in\partial\Omega$.
Such a conformal mapping will allow us to reduce the boundary blow-up analysis to
the case of a half-ball.
We derive accurate estimates for the Green's function $\mathcal G$ for the half-ball.
In Section \ref{sec:exclusion} we exclude the existence of blow-up points on $\partial\Omega$.
To this end, we argue by contradiction and we assume that $x_0\in\partial\Omega$
is a blow-up point. By extending the Brezis-Merle arguments
to the boundary via reflection arguments, we prove that a minimal mass is necessary for a boundary blow-up, so that
$x_0$ is isolated. Consequently, we use the conformal mapping $X_0$ to 
pull-back the problem to the half-ball. Exploiting the estimates for $\mathcal G$
we estimate the blow-up sequence in a small ball near the blow-up point $x_0$.
Then, a Pohozhaev identity yields the desired contradiction.
Section \ref{sec:rv} is devoted to establishing the residual vanishing property.
Extending Brezis-Merle type arguments, we show that if 
$m_+(x_0)-m_-(x_0)>4\pi\beta_+$ or $m_-(x_0)-m_+(x_0)>4\pi\beta_-$,
then $\iint_{[-1,1]\times\Omega}e^{\alpha v_k}\,\mathcal P(d\alpha) dx\to\infty$
and consequently residual vanishing holds.
On the other hand, we check that 
if $x_0\in\mathcal S_+\cap\mathcal S_-$, then the mass relation \eqref{eqn:mr}
implies $m_+(x_0)-m_-(x_0)>4\pi\beta_+$ or $m_-(x_0)-m_+(x_0)>4\pi\beta_-$,
thus concluding the proof of residual vanishing.
In Section \ref{sec:location} we prove Theorem \ref{thm:lbp}.
To this end, we comply the complex analysis argument developed by \cite{ye}.
Finally, in Section \ref{sec:mfd} we derive the corresponding results
for the problem defined on a compact Riemannian surface
without boundary.
\paragraph{Notation}
Henceforth, we omit the notation $dx$ and ${\cal P}(d\alpha)$ when it is clear from the context, 
and we do not distinguish the sequences appearing below from their subsequences. 
For the sake of simplicity, in what follows we denote $I=[-1,1]$, $I_+=[0,1]$, $I_-=[-1,0)$.

\section{Preliminaries}\label{sec:pre}

Throughout this section, we use complex notations 
by identifying $x=(x_1,x_2)\in{\bf R}^2$ with $z=x_1+\imath x_2\in{\bf C}$ in the usual way, 
where $\imath$ denotes the imaginary unit. 

Fix $x_0\in\partial\Omega$ and take $R_0>0$ such that $B_{2R_0}(x_0)\cap\partial\Omega$ is connected. 
We may assume that $x_0=0$. 
To study the problem near $x_0=0$, 
we fix $z_0\in B_{R_0}\cap\Omega$ and take a conformal mapping $X_0:B_{R_0}\cap\bar{\Omega}\rightarrow \bar{{\bf R}}_+^2$, 
where 
\[
{\bf R}_+^2=\{(X_1,X_2) \ | \ X_2>0\}, 
\]
such that
\begin{equation}\label{eqn:cmap-1}
\left\{
\begin{array}{l}
X_0(0)=0,\quad X_0(z_0)\not\in B_3^+, \\
X_0(B_{R_0}\cap\Omega)\subset {\bf R}_+^2, \\ 
X_0(B_{R_0}\cap\Omega)\supset B_2\cap{\bf R}_+^2, \\
X_0(B_{R_0}\cap\partial\Omega)\subset\partial{\bf R}_+^2,\quad X_0(B_{R_0}\cap\partial\Omega)\supset (-3,3)\times\{0\}, 
\end{array}
\right.
\end{equation}
and that
\begin{equation}\label{eqn:cmap-2}
A_0\in C^1(\bar{B}_1^+),\quad A_0\geq\delta_0
\end{equation}
for some $\delta_0>0$, where
\[
B_r^+=B_r\cap{\bf R}_+^2,\quad A_0(X)=|g'(X)|^2,\quad g=X^{-1},\quad X=X_0
\]
for $r>0$. 
This is possible if $0<R_0\ll 1$, namely, the following lemma holds. 

\begin{lem}\label{lem:cmap}
If $0<R_0\ll 1$ then there exists a conformal mapping $X_0:B_{R_0}\cap\bar{\Omega}\rightarrow \bar{{\bf R}}_+^2$ satisfying \eqref{eqn:cmap-1}-\eqref{eqn:cmap-2}. 
\end{lem}

{\it Proof}. 
By the Carath\'{e}odory theorem, there exists $w_1:\overline{B_{R_0}\cap\Omega}\rightarrow\bar{B}_1$ such that 
$w_1(z_0)=0$, $w_1(B_{R_0}\cap\Omega)=B_1$, $w_1(\partial(B_{R_0}\cap\Omega))=\partial B_1$, 
it is holomorphic in $B_{R_0}\cap\Omega$ and is homeomorphic on $\overline{B_{R_0}\cap\Omega}$.  
We may assume that 
\[
w_1(B_{R_0}\cap\partial\Omega)\ni -1,\quad w_1(B_{R_0}\cap\partial\Omega)\not\ni 1
\]
by taking a suitable rotation. 
Let 
\[
w_2(z)=-\imath\frac{z+1}{z-1},\quad w_3(z)=z-w_2\circ w_1(0). 
\]
Then we find that $X_0=L_0(w_3\circ w_2\circ w_1)$ satisfies \eqref{eqn:cmap-1} for $L_0\gg 1$. 


Since $w_1$ is injective and $|dw_1/dz|>0$ in $B_{R_0}\cap\Omega$, 
there exists a function $H=H(z)$, 
which is holomorphic in $B_{R_0}\cap\Omega$ and continuous on $\overline{B_{R_0}\cap\Omega}$, such that 
\[
w_1(z)=(z-z_0)H(z),\quad H(\zeta)\neq 0
\]
for $z\in B_{R_0}\cap\Omega$ and $\zeta\in\overline{B_{R_0}\cap\Omega}$. 
Note that there exists $0<R_0'<R_0$ such that $z_0\not\in B_{2R_0'}$
and $\log H$ is defined as a single-valued function on $\overline{B_{R_0'}\cap\Omega}$. 
For such an $R_0'$, real-valued function 
\[
U(x)=\log|z-z_0|+\log|H(z)|
\]
is the solution of 
\[
-\Delta U=0\ \mbox{in $B_{R_0'}\cap\Omega$},\quad 
U<0\ \mbox{on $\partial B_{R_0'}\cap\Omega$},\quad 
U=0\ \mbox{on $B_{R_0'}\cap\partial\Omega$}. 
\]
Since 
\begin{equation}\label{eqn:cmap-pf-1}
U\in C_{loc}^2(\overline{B_{R_0'}\cap\Omega}\setminus (\partial B_{R_0'}\cap\partial\Omega))
\end{equation}
by the elliptic regularities, and since 
$|\nabla U|>0$ on $B_{R_0'}\cap\partial\Omega$ by the Hopf lemma, 
it holds that 
\[
\left|\frac{dw_1}{dz}\right|=|w_1||\nabla U|=|\nabla U|>0\quad \mbox{on $B_{R_0'}\cap\partial\Omega$}. 
\]
Noting the conformality of $w_1$, we conclude that 
\begin{equation}\label{eqn:cmap-pf-2}
\left|\frac{dw_1}{dz}\right|>0\quad \mbox{in $\overline{B_{R_0'}\cap\Omega}\setminus (\partial B_{R_0'}\cap\partial\Omega)$}. 
\end{equation}
Consequently, \eqref{eqn:cmap-pf-1}-\eqref{eqn:cmap-pf-2} and $X_0=L_0(w_3\circ w_2\circ w_1)$ imply that 
$X_0$ satisfies \eqref{eqn:cmap-2} and $X_0\in C^1(B_{R_0'}\cap\partial\Omega)$ for $L_0\gg 1$. 

Finally, we retake $R_0$ by $R_0=R_0'/2$ if needs, and obtain the desired conformal mapping $X_0$ for $L_0\gg 1$. \hfill\qed\ \\

Now, fix $0<R_0\ll 1$ and let $(\lambda,v)$ be a solution of \eqref{eqn:neri}. 
For every function $\varphi$ defined on $B_{R_0}\cap\bar{\Omega}$, we put 
\[
\hat{\varphi}(X)=\varphi\circ g(X),\quad g=X^{-1},\quad X=X_0. 
\]
Then we have 
\[
\left\{
\begin{array}{ll}
\displaystyle 
-\Delta_X \hat{v}= cA_0(X)\int_I \alpha e^{\alpha\hat{v}} {\cal P}(d\alpha) & \mbox{in $B_{2}^+$} \\
\displaystyle 
\hat{v}=0 & \mbox{on $B_{2}\cap\partial{\bf R}_+^2$}, 
\end{array}
\right.
\]
where
\[
c=\frac{\lambda}{\iint_{I\times\Omega}e^{\alpha' v}}. 
\]\ \\

In the proof of Theorem \ref{thm:main}, we use the Brezis-Merle inequality (\cite{bm91}). 
\begin{lem}\label{lem:bm-ineq}
Let $D\subset{\bf R}^2$ be a bounded domain, and let $u=u(x)$ be a solution of 
\[
-\Delta u=f\ \mbox{in $D$},\quad u=0\ \mbox{on $\partial D$}
\]
with $f\in L^1(D)$. 
Then, for $\delta\in (0,4\pi)$, we have 
\[
\int_D \exp\left(\frac{(4\pi-\delta)|u(x)|}{\|f\|_{L^1(D)}}\right)\leq \frac{4\pi^2}{\delta}(\diam(D))^2, 
\]
where $\diam(D)=\sup_{x,y\in D}|x-y|$.  
\end{lem}\ \\

Finally, we provide the estimates concerning the Green function ${\cal G}={\cal G}(x,y)$ for $B_1^+$ defined by 
\begin{equation}\label{eqn:green-halfdisk}
{\cal G}(x,y)=-\frac{1}{2\pi}\log\left|\frac{(z-w)(1-zw)}{(z-\bar{w})(1-z\bar{w})}\right|, 
\end{equation}
that is, 
\[
-\Delta {\cal G}(\cdot,y)=\delta_y\quad \mbox{in $B_1^+$},\quad {\cal G}(\cdot,y)=0\quad \mbox{on $\partial B_1^+$}, 
\]
where $z=x_1+\imath x_2$ and $w=y_1+\imath y_2$. 

\begin{lem}\label{lem:green-esti}
Given $0<\delta\ll 1$, we have 
\begin{align}
|{\cal G}(x,y)|\leq \frac{1}{2\pi}\log\frac{2(1+\delta)}{1-\delta}\quad
\mbox{for $(x,y)\in (B_1^+\setminus B_{3\delta}^+)\times B_\delta^+$},
 \label{eqn:green-1}\\ 
|\nabla_x {\cal G}(x,y)|\leq\frac{2}{\pi(1-\delta)^2}\quad
\mbox{for $(x,y)\in (B_1^+\setminus B_{\sqrt{\delta}+\delta})\times B_\delta^+$}.
 \label{eqn:green-2}
\end{align}
\end{lem}

{\it Proof}. In this proof, we again use $z=x_1+\imath x_2$ and $w=y_1+\imath y_2$ in complex notations. 
Note that $|\nabla_x {\cal G}(x,y)|=|\frac{d}{dz}{\cal H}(z,w)|$ and 
\[
\frac{d}{dz}{\cal H}(z,w)=-\frac{1-|w|^2}{2\pi}\cdot\frac{(w-\bar{w})(1-z^2)}{(z-w)(z-\bar{w})(1-zw)(1-z\bar{w})}, 
\]
where 
\[
{\cal H}(z,w)=-\frac{1}{2\pi}\log\frac{(z-w)(1-zw)}{(z-\bar{w})(1-z\bar{w})}. 
\]

[Proof of \eqref{eqn:green-1}] 
For $(z,w)\in (B_1^+\setminus B_{3\delta}^+)\times B_\delta^+$, we compute 
\begin{align*}
\left|\frac{(z-w)(1-zw)}{(z-\bar{w})(1-z\bar{w})}\right|
&\leq\frac{|z|+|w|}{|z|-|w|}\cdot\frac{1+|w||z|}{1-|w||z|}
\leq\frac{|z|+\delta}{|z|-\delta}\cdot\frac{1+\delta|z|}{1-\delta|z|}\\
&\leq\frac{1+\delta}{1-\delta}\left(1+\frac{2\delta}{|z|-\delta}\right)
\leq\frac{2(1+\delta)}{1-\delta}
\end{align*}
Similarly, for $(z,w)\in (B_1^+\setminus B_{3\delta}^+)\times B_\delta^+$, we compute
\[
\left|\frac{(z-w)(1-zw)}{(z-\bar{w})(1-z\bar{w})}\right|\geq \left(\frac{2(1+\delta)}{1-\delta}\right)^{-1}. 
\]

[Proof of \eqref{eqn:green-2}] 
For $(z,w)\in (B_1^+\setminus B_{\sqrt{\delta}+\delta}^+)\times B_\delta^+$, we compute 
\begin{align*}
\left|\frac{d{\cal H}}{dz}(z,w)\right|
&\leq\frac{1}{2\pi}\cdot\frac{4|w|}{(|z|-|w|)^2(1-|z||w|)^2}\\
&\leq\frac{2\delta}{\pi(|z|-\delta)^2(1-\delta|z|)^2}
\leq\frac{2}{\pi(1-\delta)^2}. 
\end{align*}

Hence, the desired estimates are shown. \hfill\qed\
\section{Proof of Theorem \ref{thm:main}: proof of \eqref{eqn:exclusion}}\label{sec:exclusion}

To begin with, we collect in Proposition \ref{prop:rz2012} below 
some results from \cite{rz2012} which are used in the proof of Theorem \ref{thm:main}.
Although such results were actually derived in the case of a compact manifold, 
the extension to the case of \textit{interior} blow-up points for the Dirichlet problem is straightforward
since the arguments are local in nature.
\par
We introduce the measure functions $\nu_{k,\pm}=\nu_{k,\pm}(dx)\in{\cal M}(\bar{\Omega})$ 
and $\mu_k=\mu_k(d\alpha dx)\in{\cal M}(I\times\bar{\Omega})$ defined by 
\begin{align}
&\nu_{k,\pm}=\lk\int_{I_\pm}|\alpha| V(\alpha,v_k)e^{\alpha \vk}{\cal P}(d\alpha), 
 \label{eqn:def-nuk}\\
&\mu_k(d\alpha dx)=\lk V(\alpha,v_k)e^{\alpha \vk}{\cal P}(d\alpha)dx, 
\end{align}
respectively, where $V(\alpha,v)=\left(\iint_{[-1,1]\times\Omega}e^{\alpha v}\,\mathcal P(d\alpha)dx\right)^{-1}$
if $\vk$ satisfies \eqref{eqn:neri} and $V(\alpha,v)=\left(\int_\Omega e^{\alpha v}\,dx\right)^{-1}$
if $\vk$ satisfies \eqref{eqn:ss}.
We recall that $I_+=[0,1]$, $I_-=[-1,0)$. 
Then, Theorem~2.1 and Theorem~2.2 in \cite{rz2012} are readily adapted to the domain case to yield the following
basic blow-up properties. 
\begin{pro}
[\cite{rz2012}, common blow-up properties for \eqref{eqn:neri} and \eqref{eqn:ss}]
\label{prop:rz2012}
\label{fact:rz}
Let $(\lk,\vk)$ be a solution sequence of \eqref{eqn:neri} or \eqref{eqn:ss} with $\lk\rightarrow\lzero$ for some $\lzero\geq 0$. 
Assume that 
\[
\mathcal S\cap\partial\Omega=\emptyset
\]
holds true. 
Then, passing to a subsequence, we have the following alternatives. 
\begin{enumerate}
\item[{\bf (I)}] Compactness: $\displaystyle \limsup_{k\rightarrow\infty}\|\vk\|_\infty<+\infty$, that is, $\cS=\emptyset$. \\
Then, there exists $v\in H_0^1(\Omega)$ such that $\vk\rightarrow v$ in $H_0^1(\Omega)$ and $v$ is a solution of \eqref{eqn:neri}
or \eqref{eqn:ss}. 
\item[{\bf (II)}] Concentration: $\displaystyle \limsup_{k\rightarrow\infty}\|\vk\|_\infty=\infty$, that is, $\cS\neq\emptyset$. \\
Then, $\cS$ is finite and the following properties {\bf a)}-{\bf c)} hold:
\begin{enumerate}
\item[{\bf a)}] There exists $0\leq s_\pm\in L^1(\Omega)\cap L_{loc}^\infty(\Omega\setminus\cS)$ such that
\begin{equation}\label{eqn:rz1}
\nu_{k,\pm} \overset{\ast}{\rightharpoonup} \nu_\pm=s_\pm+\sum_{x_0\in\cS_\pm}m_\pm(x_0)\delta_{x_0} \quad \mbox{in ${\cal M}(\bar{\Omega})$}, 
\end{equation}
with $m_\pm(x_0)\geq 4\pi$ for every $x_0\in\cS_\pm$, 
where $\delta_{x}\in{\cal M}(\Omega)$, $x\in\Omega$, denotes the Dirac measure centered at $x$. 
\item[{\bf b)}] There exist $\zeta_{x_0}\in{\cal M}(I)$ and $0\leq r\in L^1(I\times\Omega)$ such that 
\begin{equation}\label{eqn:rz2}
\mu_k \overset{\ast}{\rightharpoonup} 
\mu=\mu(d\alpha dx)
=r(\alpha,x){\cal P}(d\alpha)dx+\sum_{x_0\in\cS}\zeta_{x_0}(d\alpha)\delta_{x_0}(dx) \quad \mbox{in ${\cal M}(I\times\bar{\Omega})$}. 
\end{equation}
\item[{\bf c)}] For every $x_0\in\cS$, we have  
\begin{align}
&8\pi\int_I \zeta_{x_0}(d\alpha)=\left(\int_I \alpha \zeta_{x_0}(d\alpha)\right)^2
 \label{eqn:rz3}\\
&m_\pm(x_0)=\int_{I_\pm}|\alpha| \zeta_{x_0}(d\alpha),\quad s_\pm(x)=\int_{I_\pm}|\alpha| r(\alpha,x){\cal P}(d\alpha), 
 \label{eqn:rz4}
\end{align}
where $m_\pm(x_0)$ is as in \eqref{eqn:rz1}. 
Moreover, for every $x_0 \in\cS_\pm\setminus\cS_\mp$, it holds that
\begin{equation}\label{eqn:rz5}
m_\mp(x_0)=\int_{I_\mp}|\alpha|\zeta_{x_0}(d\alpha)=0. 
\end{equation}
\end{enumerate}
\end{enumerate}
\end{pro}\ \\

To prove Theorem \ref{thm:main}, we see from Proposition \ref{prop:rz2012} that it suffices to show 
\eqref{eqn:exclusion}, \eqref{eqn:residual-vanishing} and \eqref{eqn:mr} under the assumption that $\cS\neq\emptyset$.

In the remainder of this section, we shall prove \eqref{eqn:exclusion} by contradiction 
under the assumption that $\cS\cap\partial\Omega\neq\emptyset$. \\ 

Let $(\lk,\vk)$ be a solution sequence of \eqref{eqn:neri}. 
The starting point for the proof of \eqref{eqn:exclusion} is the following lemma 
based on the Brezis-Merle inequality, see Lemma \ref{lem:bm-ineq}.

\begin{lem}\label{lem:assurance-finiteness}
For any $x_0\in\cSpm\cap\partial\Omega$, it holds that
\begin{equation}\label{eqn:massesti-below}
\lim_{r\downarrow 0}\limsup_{k\rightarrow\infty} \nu_{k,\pm}(B_r(x_0)\cap\Omega)\geq 4\pi, 
\end{equation}
where $\nu_{k,\pm}$ is as in \eqref{eqn:def-nuk}. 
\end{lem}

{\it Proof}. Fix $x_0\in\cS\cap\partial\Omega$. 
We may assume that $x_0=0$. 
Assume that \eqref{eqn:massesti-below} is false. 
Then, there exist $0<\varepsilon_0, r_0\ll 1$ such that $B_{3r_0}\cap(\cS\cap\Omega)=\emptyset$ and 
\begin{equation}\label{eqn:lem:assurance-finiteness-1}
\nu_{k,\pm}(B_{3r_0}\cap\Omega)\leq 4\pi-2\varepsilon_0
\end{equation}
for $k\gg 1$. 
We decompose $\vk$ as $\vk=\vk^+-\vk^-$, where $\vk^\pm$ is the solution of 
\[
-\Delta \vk^\pm=\nu_{k,\pm}\ \mbox{in $\Omega$},\quad \vk^\pm=0\ \mbox{on $\partial\Omega$}. 
\]
Let $v_{k,1}^\pm$ and $v_{k,2}^\pm$ be the solutions of 
\begin{align*}
&-\Delta v_{k,1}^\pm=\nu_{k,\pm}\ \mbox{in $B_{2r_0}\cap\Omega$},\quad 
v_{k,1}^\pm=0\ \mbox{on $\partial(B_{2r_0}\cap\Omega)$}, \\
&-\Delta v_{k,2}^\pm=0\ \mbox{in $B_{2r_0}\cap\Omega$},\quad 
v_{k,2}^\pm=\vk^\pm\ \mbox{on $\partial(B_{2r_0}\cap\Omega)$},
\end{align*}
respectively. 
Then it holds that $\vk^\pm=v_{k,1}^\pm+v_{k,2}^\pm$ on $\bar{B}_{2r_0}\cap\bar{\Omega}$. 
In addition, 
by the maximum principle and the $L^1$-estimates (see \cite{bs73}),  
we have 
\begin{align}
&\vk^\pm\geq 0\ \mbox{on $\bar{\Omega}$},\quad v_{k,1}^\pm\geq 0\ \mbox{on $\bar{B}_{2r_0}\cap\bar{\Omega}$}, 
 \label{eqn:lem:assurance-finiteness-2}\\
&\|\vk^\pm\|_{L^1(\Omega)}+\|v_{k,1}^\pm\|_{L^1(B_{2r_0}\cap\Omega)}\leq C_{\const\insz} 
 \label{eqn:lem:assurance-finiteness-3}
\end{align}
for some $C_{\subz}>0$ independent of $k$ and $0<r_0\ll 1$. 

To estimate $v_{k,2}^\pm=v_{k,2}^\pm(x)$, 
we take a conformal mapping $X:B_{2r_0}\cap\bar{\Omega}\rightarrow\bar{{\bf R}}^2$ as in Section \ref{sec:pre} for $R_0=2r_0$. 
Let $\tilde{v}_{k,2}^\pm$ be the odd extension of $v_{k,2}^\pm\circ X^{-1}$. 
Then it holds that 
\[
-\Delta_X \tilde{v}_{k,2}^\pm=0 \ \mbox{in $B_1$},\quad 
\tilde{v}_{k,2}^\pm=0 \ \mbox{on $B_1\cap\partial{\bf R}_+^2$}, 
\]
and the mean value theorem and \eqref{eqn:lem:assurance-finiteness-3} admit $C_{\const\insy}>0$, independent of $k$, such that  
\begin{equation}\label{eqn:lem:assurance-finiteness-4}
\|v_{k,2}^\pm\|_{L^\infty(B_{r_1}\cap\Omega)}\leq C_{\suby}
\end{equation}
for $k\gg 1$ and for some $0<r_1<2r_0$. 

On the other hand, Lemma \ref{lem:bm-ineq} and \eqref{eqn:lem:assurance-finiteness-1} yield 
\begin{equation}\label{eqn:lem:assurance-finiteness-5}
\int_{B_{2r_0}\cap\Omega} \exp\left(\frac{4\pi-\varepsilon_0}{4\pi-2\varepsilon_0}v_{k,1}^\pm\right)
\leq \frac{4\pi^2}{\varepsilon_0}(4r_0)^2
\end{equation}
for $k\gg 1$. 
Combining \eqref{eqn:lem:assurance-finiteness-2} and  \eqref{eqn:lem:assurance-finiteness-4}-\eqref{eqn:lem:assurance-finiteness-5}, and noting that $\pm\vk\leq \vk^\pm$, we obtain 
\begin{align*}
&\int_{B_{r_1}\cap\Omega} e^{p_0 \alpha \vk}
\leq\int_{B_{r_1}\cap\Omega} e^{p_0 \vk^\pm} \\
&\leq \int_{B_{r_1}\cap\Omega} e^{p_0(v_{k,1}^\pm+|v_{k,2}^\pm|)}
\leq \frac{(8\pi r_0)^2}{\varepsilon_0}\cdot e^{p_0 C_{\suby}}
\end{align*}
for any $\alpha\in I_\pm$, where $p_0=(4\pi-\varepsilon_0)/(4\pi-2\varepsilon_0)>1$. 
This estimate means 
\[
\|\nu_{k,\pm}\|_{L^{p_0}(B_{r_1}\cap\Omega)}^{p_0} \leq C_{\const\insx}  
\]
for any $k\gg 1$ and for some $C_{\subx}>0$ independent of $k$. 
Consequently, the boundary $L^p$-estimate guarantees the uniform boundedness of $\vk^\pm$ in $B_{r_1/2}\cap\Omega$, 
which contradicts $0\in\cSpm$. \hfill\qed\ \\

Now we fix $x_0\in\cS\cap\partial\Omega$. 
We may assume that $x_0=0$. 
By virtue of Lemma \ref{lem:assurance-finiteness}, we see that $\cS$ is finite, 
and hence there exists $0<R_0\ll 1$ such that $B_{2R_0}\cap \cS=\{0\}$. 
After passing to a subsequence, we set 
\[
m(0)=\lim_{r\downarrow 0}\lim_{k\rightarrow\infty}
\nu_{k,+}(B_r\cap\Omega)+\nu_{k,-}(B_r\cap\Omega)\geq 4\pi. 
\]

Given $0<\varepsilon\ll 1$, there exists $r_\varepsilon\in(0,2R_0)$ such that 
\begin{equation*}
\lim_{k\rightarrow\infty}
\nu_{k,+}(B_{r_\varepsilon}\cap\Omega)+\nu_{k,-}(B_{r_\varepsilon}\cap\Omega)
\leq m(0)+\varepsilon/4. 
\end{equation*}
We transform the problem into the one on $B_1^+$ 
by taking the conformal mapping $X:B_{r_\varepsilon}\cap\bar{\Omega}\rightarrow \bar{{\bf R}}_+^2$ as in Section \ref{sec:pre} for $R_0=r_\varepsilon$. 
For simplicity, we shall denote $\hat{v}_k$, $X$, $\nabla_X$ and $\Delta_X$ 
by $\vk$, $x$, $\nabla$ and $\Delta$, respectively. 
Under these agreements, we obtain 
\begin{equation}\label{eqn:neri-semi-disk}
\left\{
\begin{array}{ll}
-\Delta \vk=\nu_k & \mbox{in $B_{2}^+$} \\
\displaystyle 
\vk=0 & \mbox{on $B_{2}\cap\partial{\bf R}_+^2$}, 
\end{array}
\right.
\end{equation}
where 
\[
\nu_k=\nu_{k,+}-\nu_{k,-}=c_kA_0(x)\int_I\alpha e^{\alpha\vk}{\cal P}(d\alpha),\quad  
c_k=\frac{\lk}{\iint_{I\times\Omega}e^{\alpha' \vk}}. 
\]
Note that 
\begin{align}
&m(0)=\lim_{r\downarrow 0}\lim_{k\rightarrow\infty}\nu_{k,+}(B_r^+)+\nu_{k,-}(B_r^+)\geq 4\pi, \label{eqn:mass-1} \\
&\lim_{k\rightarrow\infty}\nu_{k,+}(B_1^+)+\nu_{k,-}(B_1^+)\leq m(0)+\varepsilon/4, 
 \label{eqn:mass-2}
\end{align}
and that there exists $r_k\downarrow 0$ such that 
\begin{equation}\label{eqn:mass-3}
\lim_{k\rightarrow\infty}\nu_{k,+}(B_{r_k}^+)+\nu_{k,-}(B_{r_k}^+)=m(0). 
\end{equation}

We now show the crucial estimate. 

\begin{lem}\label{lem:crucial-estimate}
There exists $C_{\const\insa}>0$, independent of $\delta$, such that  
\[
\limsup_{k\rightarrow\infty}\|v_k\|_{W^{1,\infty}(B_1^+\setminus B_\delta^+)}\leq C_{\suba}
\]
for any $0<\delta\ll 1$. 
\end{lem}

{\it Proof}. The proof is split into five steps. \\

{\it Step 1}. Let $\wk$ and $\hk$ be the solutions of 
\[
\left\{
\begin{aligned}
-\Delta\wk=&\nu_k &&\mbox{in $B_1^+$} \\
\wk=&0 &&\mbox{on $\partial B_1^+$},
\end{aligned}
\right.\quad 
\left\{
\begin{aligned}
-\Delta\hk=&0 &&\mbox{in $B_1^+$} \\
\hk=&\vk &&\mbox{on $\partial B_1^+$}, 
\end{aligned}
\right.
\]
respectively, so that $\vk=\wk+\hk$. 
Then, there exists $C_{\const\insb}>0$, independent of $k$, such that  
\begin{equation}\label{eqn:crucial-1}
\|\hk\|_{W^{1,\infty}(B_1^+)}\leq C_{\subb}. 
\end{equation}
In fact, the odd extension $\thk=\thk(x)$ of $\hk$ is the solution of
\[
-\Delta\thk=0\ \mbox{in $B_1$},\quad \thk=\tvk\ \mbox{on $\partial B_1$}, 
\]
where $\tvk=\tvk(x)$ is the odd one of $\vk$. 
Therefore, $\|\thk\|_{W^{1,\infty}(B_1)}$ is uniformly bounded by the maximum principle 
and the uniform boundedness of $\|\vk\|_{W^{1,\infty}(\partial B_1^+)}$, which means \eqref{eqn:crucial-1}. \\

{\it Step 2}. Let $\wkn$ and $\wknn$ be the solutions of 
\[
\left\{
\begin{aligned}
-\Delta\wkn=&\chi_{B_{\rk}^+}\nu_{k} &&\mbox{in $B_1^+$} \\
\wkn=&0 &&\mbox{on $\partial B_1^+$},
\end{aligned}
\right.\quad
\left\{
\begin{aligned}
-\Delta\wknn=&\chi_{B_1^+\setminus B_{\rk}^+}\nu_{k} &&\mbox{in $B_1^+$} \\
\wknn=&0 &&\mbox{on $\partial B_1^+$},
\end{aligned}
\right.
\]
so that $\wk=\wkn+\wknn$, recall that $\rk$ satisfies \eqref{eqn:mass-3}, 
where $\chi_A$ denotes the characteristic function of $A\subset\Rd$. 
Let $\Gg=\Gg(x,y)$ be the Green function for $B_1^+$ defined by \eqref{eqn:green-halfdisk}. 
Then, the representation formula
\[
\wkn(x)=\int_{B_1^+}\Gg(x,y)\chi_{B_{\rk}^+}(y)\nu_{k}(y)dy
=\int_{B_{\rk}^+}\Gg(x,y)\nu_{k}(y)dy
\]
shows 
\begin{equation}\label{eqn:modi-exclusion-1}
\lim_{k\rightarrow\infty}\wkn(x)=0
\end{equation}
for every $x\in \bar{B}_1^+$, 
since $\rk\downarrow 0$, 
$\lim_{y\rightarrow\partial B_1^+}\Gg(x,y)=0$ and $\Gg(x,\cdot)\in C(\bar{\Omega}\setminus\{x\})$ for $x\in B_1^+$, 
and $\|\nu_k\|_{L^1(B_1^+)}$ is uniformly bounded. 

Fix $0<\varepsilon\ll 1$. 
Then, the representation formula of $\wkn$ above and \eqref{eqn:green-1} admit $C_{\const\insp}>0$, 
independent of $\varepsilon$, 
such that 
\begin{equation}\label{eqn:modi-exclusion-2}
\limsup_{k\rightarrow\infty}\|\wkn\|_{L^\infty(B_1^+\setminus B_{\varepsilon/2}^+)}\leq C_{\subp}. 
\end{equation}
Here, we consider the problem 
\[
\left\{
\begin{aligned}
-\Delta\psik=&0 &&\mbox{in $B_1^+\setminus B_\varepsilon^+$} \\
\psik=&\wkn &&\mbox{on $\partial B_\varepsilon\cap{\bf R}_+^2$} \\
\psik=&0 &&\mbox{on $\partial(B_1^+\setminus B_\varepsilon^+)\setminus(\partial B_\varepsilon\cap{\bf R}_+^2)$}. 
\end{aligned}
\right.
\]
Note that $\psik=\wkn$ on $\overline{B_1^+\setminus B_\varepsilon^+}$ for $k\gg 1$ since $\rk\downarrow 0$. 
Let $\tpsik=\tpsik(x)$ and $\twkn=\twkn(x)$ be the odd extensions of $\psik$ and $\wkn$, respectively. 
Then, we have 
\[
\left\{
\begin{aligned}
-\Delta\tpsik=&0 &&\mbox{in $B_1\setminus B_\varepsilon$} \\
\tpsik=&\twkn &&\mbox{on $\partial B_\varepsilon$} \\
\tpsik=&0 &&\mbox{on $\partial B_1$}
\end{aligned}
\right.
\]
and use the maximum principle and \eqref{eqn:modi-exclusion-2} 
to find that there exists $C_{\const\insq}>0$, independent of $\varepsilon$, such that  
\[
\limsup_{k\rightarrow\infty}\|\twkn\|_{L^\infty(B_1\setminus B_\varepsilon)}
=\limsup_{k\rightarrow\infty}\|\tpsik\|_{L^\infty(B_1\setminus B_\varepsilon)}
\leq C_{\subq}, 
\]
where we have used the property that $\psik=\wkn$ on $\overline{B_1^+\setminus B_\varepsilon^+}$ for $k\gg 1$. 
Thus, the elliptic regularity yields $C_{\const\insr,\varepsilon}>0$ such that  
\[
\limsup_{k\rightarrow\infty}\|\twkn\|_{W^{1,\infty}(B_1\setminus B_{2\varepsilon})}
=\limsup_{k\rightarrow\infty}\|\tpsik\|_{W^{1,\infty}(B_1\setminus B_{2\varepsilon})}
\leq C_{\subr,\varepsilon}. 
\]
Hence, the Arzel\`{a}-Ascoli theorem and \eqref{eqn:modi-exclusion-1} guarantee that 
\[
\wkn\rightarrow 0\quad \mbox{in $C(\bar{B}_1^+\setminus B_{2\varepsilon})$}. 
\]
Since $\varepsilon$ is arbitrary, we conclude that 
\begin{equation}\label{eqn:crucial-2}
\wkn\rightarrow 0\quad \mbox{in $C_{loc}(\bar{B}_1^+\setminus\{0\})$}. 
\end{equation}

{\it Step 3}. Now \eqref{eqn:mass-1}-\eqref{eqn:mass-2} and Lemma \ref{lem:bm-ineq} show 
\begin{equation}\label{eqn:crucial-3}
\int_{B_1^+}e^{2|\wknn|}
=\int_{B_1^+}e^{\frac{4\pi-(4\pi-\varepsilon)}{\varepsilon/2}|\wknn|}
\leq \int_{B_1^+}e^{\frac{4\pi-(4\pi-\varepsilon)}{\|\chi_{B_1^+\setminus B_{\rk}^+}\nu_k\|_{L^1(B_1^+)}}|\wknn|}
\leq 8\pi
\end{equation}
for $k\gg 1$. 
Summarizing \eqref{eqn:crucial-1}, \eqref{eqn:crucial-2}-\eqref{eqn:crucial-3} and the uniform boundedness of $c_k A_0$, 
we obtain $C_{\const\insc}>0$, independent of $\delta$, such that 
\begin{equation}\label{eqn:crucial-4}
\limsup_{k\rightarrow\infty}\|\nu_k\|_{L^2(B_1^+\setminus B_\delta^+)}\leq C_{\subc}
\end{equation}
for any $0<\delta\ll 1$. \\

{\it Step 4}. In this step, we shall derive the $L^\infty$-estimates of $\vk$. 
Given $0<\delta\ll 1$, let $\zkd$ and $\zkkd$ be the solutions of 
\begin{equation}\label{eqn:zkd-zkkd}
\left\{
\begin{aligned}
-\Delta\zkd=&\chi_{B_\delta^+}\nu_k &&\mbox{in $B_1^+$} \\
\zkd=&0 &&\mbox{on $\partial B_1^+$},
\end{aligned}
\right.\quad
\left\{
\begin{aligned}
-\Delta\zkkd=&\chi_{B_1^+\setminus B_\delta^+}\nu_k &&\mbox{in $B_1^+$} \\
\zkkd=&0 &&\mbox{on $\partial B_1^+$}, 
\end{aligned}
\right.
\end{equation}
so that $\wk=\zkd+\zkkd$. 
It follows from \eqref{eqn:crucial-4} and the elliptic regularity that there exists $C_{\const\insd}>0$, independent of $\delta$, such that 
\begin{equation}\label{eqn:crucial-5}
\limsup_{k\rightarrow\infty}\|\zkkd\|_{L^\infty(B_1^+)}\leq C_{\subd}. 
\end{equation}
Furthermore, the representation formula
\[
\zkd(x)=\int_{B_1^+}\Gg(x,y)\chi_{B_\delta^+}(y)\nu_k(y)dy=\int_{B_\delta^+}\Gg(x,y)\nu_k(y)dy
\]
and \eqref{eqn:green-1} admit $C_{\const\inse}>0$, independent of $\delta$, such that 
\begin{equation}\label{eqn:crucial-6}
\limsup_{k\rightarrow\infty}\|\zkd\|_{L^\infty(B_1^+\setminus B_{3\delta}^+)}\leq C_{\sube}. 
\end{equation}
Since $0<\delta\ll 1$ is arbitrary and since the constants in \eqref{eqn:crucial-1} and \eqref{eqn:crucial-5}-\eqref{eqn:crucial-6} are independent of $\delta$, 
we conclude that there exists $C_{\const\insf}>0$, independent of $\delta'$, such that  
\begin{equation}\label{eqn:crucial-7}
\limsup_{k\rightarrow\infty}\|\vk\|_{L^\infty(B_1^+\setminus B_{\delta'}^+)}\leq C_{\subf}
\end{equation}
for any $0<\delta'\ll 1$. \\ 

{\it Step 5}.  In the final step, we shall derive the gradient estimates of $\vk$. 
Given $0<\delta\ll 1$, we again use the decomposition $\wk=\zkd+\zkkd$ or \eqref{eqn:zkd-zkkd}. 
From \eqref{eqn:crucial-7} and the uniform boundedness of $c_k A_0$, 
we see that there exists $C_{\const\insz}>0$, independent of $\delta$, such that 
\[
\limsup_{k\rightarrow\infty}\|\chi_{B_1^+\setminus B_\delta^+} \nu_k\|_{L^\infty(B_1^+)}\leq C_{\subz}. 
\]
Thus, the elliptic regularity yields $C_{\const\insg}>0$, independent of $\delta$, such that  
\begin{equation}\label{eqn:crucial-8}
\limsup_{k\rightarrow\infty}\|\zkkd\|_{W^{1,\infty}(B_1^+)}\leq C_{\subg}. 
\end{equation}
By \eqref{eqn:green-2}, \eqref{eqn:mass-2} and the representation formula
\[
\nabla\zkd(x)=\int_{B_1^+}\chi_{B_\delta^+}(y)\nu_{k}(y)\nabla_x\Gg(x,y)dy
=\int_{B_\delta^+}\nu_{k}(y)\nabla_x\Gg(x,y)dy 
\]
it holds that 
\begin{equation}\label{eqn:crucial-9}
|\nabla\zkd(x)|\leq \frac{2(m(0)+\varepsilon)}{\pi(1-\delta)^2}
\end{equation}
for any $x\in B_1^+\setminus B_{\sqrt{\delta}+\delta}^+$ and $k\gg 1$. 

Since $\delta$ is arbitrary, 
we combine \eqref{eqn:crucial-8}-\eqref{eqn:crucial-9} and \eqref{eqn:crucial-1} to obtain $C_{\const\insh}>0$,  
independent of $\delta'$, such that 
\[
\limsup_{k\rightarrow\infty}\|\nabla\vk\|_{L^\infty(B_1^+\setminus B_{\delta'}^+)}\leq C_{\subh}
\]
for any $0<\delta'\ll 1$. 
Hence, the desired gradient estimates are established. \hfill\qed\\ 

We are now in a position to prove \eqref{eqn:exclusion}. \\

\underline{Proof of \eqref{eqn:exclusion}}\ 
At first, we note the following Pohozaev type identity, that is, 
\begin{align}\label{eqn:pohozaev}
&r\int_{\partial B_r^+\cap{\bf R}^2}\left(\frac{|\nabla\vk|^2}{2}-({\bf n}\cdot\nabla\vk)^2\right)-c_kA_0e^{\alpha\vk}d\sigma \nonumber\\
&=-\int_{B_r^+}\left\{2+\left(x\cdot\frac{\nabla A_0}{A_0}\right)\right\}c_kA_0\left(\int_I e^{\alpha\vk}{\cal P}(d\alpha)\right)dx
\end{align}
for any $0<r<1$, where $\cdot$, ${\bf n}$ and $d\sigma$ denote the usual inner product in $\Rd$, the outward unit normal vector to the boundary, and the line element on the boundary, respectively. 
Identity \eqref{eqn:pohozaev} is shown by multiplying \eqref{eqn:neri-semi-disk} by $x\cdot\nabla \vk$ without difficulty. 

Next, organizing \eqref{eqn:pohozaev}, the uniform boundedness of $c_kA_0$, \eqref{eqn:cmap-2} and Lemma \ref{lem:crucial-estimate}, 
we find  
\begin{equation}\label{eqn:exclusion-proof}
O(r^2)=-(2+O(r))
\lim_{k\rightarrow\infty}\int_{B_r^+}c_k A_0\left(\int_I e^{\alpha\vk}{\cal P}(d\alpha)\right)dx
\end{equation}
as $r\downarrow 0$ after taking $k\rightarrow\infty$ and passing to a subsequence. 
However, the right-hand side of \eqref{eqn:exclusion-proof} does not converge to $0$ as $r\downarrow 0$ because of \eqref{eqn:mass-1}, which is a contradiction. 
The proof is complete. \hfill\qed

\section{Proof of Theorem \ref{thm:main}: proof of \eqref{eqn:residual-vanishing} and \eqref{eqn:mr}}\label{sec:rv}

In this section, we shall give the proof of \eqref{eqn:residual-vanishing} and \eqref{eqn:mr} 
under the assumption that the alternative {\bf (II)} in Proposition \ref{prop:rz2012} occurs. 
Note that $\cS\cap\partial\Omega=\emptyset$ as shown in the previous section. 
Therefore, \eqref{eqn:rz1}-\eqref{eqn:rz5} now hold. \\ 

At first we shall prove \eqref{eqn:mr}. 
Suppose that the following propositions hold. 
\begin{align}
&x_0\in\cS_\pm\setminus\cS_\mp \ \Rightarrow \ \supp\zeta_{x_0}=\{\alpha_\pm^\ast\},
 \label{eqn:supp1}\\
&x_0\in\cSp\cap\cSm \ \Rightarrow \ \supp\zeta_{x_0}=\{\alpha_+^\ast,\alpha_-^\ast\}, 
 \label{eqn:supp2}
\end{align}
where $\alpha_\pm^\ast$ is as in \eqref{def:suppP}. 
Note that $\alpha_\pm^\ast\neq 0$ if $\cS_\pm\neq\emptyset$. 
Then, \eqref{eqn:rz3}-\eqref{eqn:rz5} imply that 
\[
m_\pm(x_0)=8\pi\beta_\pm,\quad m_\mp(x_0)=0
\]
if $x_0\in\cS_\pm\setminus\cS_\mp$, and that 
\[
(m_+(x_0)-m_-(x_0))^2=8\pi(\beta_+ m_+(x_0)+\beta_- m_-(x_0))
\]
if $x_0\in\cSp\cap\cSm$, where $\beta_\pm$ is as in \eqref{eqn:mq4}. 
Therefore \eqref{eqn:supp1}-\eqref{eqn:supp2} assure \eqref{eqn:mr}, 
and hence the proof of \eqref{eqn:mr} is reduced to showing \eqref{eqn:supp1}-\eqref{eqn:supp2}. \\

\underline{{\it Proof of \eqref{eqn:mr}}} 
It suffices to prove \eqref{eqn:supp1}-\eqref{eqn:supp2} as we have seen above. 
We shall only give the proof of \eqref{eqn:supp2} here since that of \eqref{eqn:supp1} is similar. 

We now fix $x_0\in\cSp\cap\cSm$. 
Then it holds that $\alpha_+^\ast>0>\alpha_-^\ast$. 
Furthermore, we find that the proof of \eqref{eqn:supp2} is reduced to proving that 
for any $0<\varepsilon\ll 1$, there exists $C_{\const\insa}>0$ such that 
\begin{equation}\label{eqn:mq-1}
\left\|\frac{e^{\alpha \vk}}{\iint_{I\times\Omega}e^{\alpha' \vk}}\right\|_{L^{p_\varepsilon}(\Omega)}\leq C_{\suba}
\end{equation}
for any $k$ and $\alpha\in [\alpha_-^\ast+2\varepsilon,\alpha_+^\ast-2\varepsilon]$, where 
\[
p_\varepsilon=\min\left\{\frac{\alpha_+^\ast-\varepsilon}{\alpha_+^\ast-2\varepsilon},\frac{\alpha_-^\ast+\varepsilon}{\alpha_-^\ast+2\varepsilon}\right\}>1. 
\]
In fact, if \eqref{eqn:mq-1} holds 
then $\zeta_{x_0}([\alpha_-^\ast+2\varepsilon,\alpha_+^\ast-2\varepsilon])=0$ for any $0<\varepsilon\ll 1$, which implies \eqref{eqn:supp2}. 
Inequality \eqref{eqn:mq-1} is obvious for $\alpha=0$. 

Given $\alpha\in (0,\alpha_+^\ast-2\varepsilon]$, we have 
\begin{equation}\label{eqn:mq-2}
\int_\Omega e^{p_\varepsilon \alpha\vk}
\leq |\Omega|^{1-\frac{\alpha p_\varepsilon}{\beta}}\left(\int_\Omega e^{\beta \vk}\right)^{\frac{\alpha p_\varepsilon}{\beta}}
\leq |\Omega|+\int_\Omega e^{\beta \vk}
\end{equation}
for any $\beta\in [\alpha_+^\ast-\varepsilon,\alpha_+^\ast]$. 
Note that $\beta/(\alpha p_\varepsilon)\geq 1$ for any $\beta\in [\alpha_+^\ast-\varepsilon,\alpha_+^\ast]$. 
Then, we combine \eqref{eqn:mq-2} with ${\cal P}([\alpha_+^\ast-\varepsilon,\alpha_+^\ast])>0$ to obtain
\[
\int_\Omega e^{p_\varepsilon \alpha\vk}\leq |\Omega|+\frac{1}{{\cal P}([\alpha_+^\ast-\varepsilon,\alpha_+^\ast])}
\iint_{I\times\Omega} e^{\alpha' \vk} {\cal P}(d\alpha')dx, 
\]
that is, 
\begin{equation}\label{eqn:mq-3}
\left\|\frac{e^{\alpha \vk}}{\iint_{I\times\Omega}e^{\alpha' \vk}}\right\|_{L^{p_\varepsilon}(\Omega)}^{p_\varepsilon}
\leq \frac{|\Omega|}{\left(\iint_{I\times\Omega}e^{\alpha' \vk}\right)^{p_\varepsilon}}
+\frac{1}{{\cal P}([\alpha_+^\ast-\varepsilon,\alpha_+^\ast])\left(\iint_{I\times\Omega}e^{\alpha' \vk}\right)^{p_\varepsilon-1}}
\end{equation}
for any $\alpha\in (0,\alpha_+^\ast-2\varepsilon]$. 
The similar argument yields 
\begin{equation}\label{eqn:mq-3'}
\left\|\frac{e^{\alpha \vk}}{\iint_{I\times\Omega}e^{\alpha' \vk}}\right\|_{L^{p_\varepsilon}(\Omega)}^{p_\varepsilon}
\leq \frac{|\Omega|}{\left(\iint_{I\times\Omega}e^{\alpha' \vk}\right)^{p_\varepsilon}}
+\frac{1}{{\cal P}([\alpha_-^\ast,\alpha_-^\ast+\varepsilon])\left(\iint_{I\times\Omega}e^{\alpha' \vk}\right)^{p_\varepsilon-1}}
\end{equation}
for any $\alpha\in [\alpha_-^\ast+2\varepsilon,0)$. 

On the other hand, for $\omega\subset\subset\Omega\setminus\cS$, there exists $C_{\const\insb,\omega}>0$ such that 
\[
\|\vk\|_{L^\infty(\omega)}\leq C_{\subb,\omega}
\]
for any $k$, and thus 
\begin{equation}\label{eqn:mq-4}
\iint_{I\times\Omega}e^{\alpha' \vk}\geq \iint_{I\times\omega}e^{\alpha' \vk}\geq |\omega|e^{-C_{\subb,\omega}}>0
\end{equation}
for any $k$. 
Consequently, \eqref{eqn:mq-1} follows from \eqref{eqn:mq-3}-\eqref{eqn:mq-4} \qed \ \\

It is left to prove \eqref{eqn:residual-vanishing}. 
For the purpose, we prepare the following lemma. 

\begin{lem}\label{lem:judge-rv}
For $x_0\in\cS$, if 
\begin{equation}\label{eqn:pf-rv-1}
m_+(x_0)-m_-(x_0)>4\pi\beta_+ \quad\mbox{or}\quad m_-(x_0)-m_+(x_0)>4\pi\beta_-,
\end{equation}
then 
\begin{equation}\label{eqn:lim-rv}
\lim_{k\rightarrow\infty}\iint_{I\times\Omega}e^{\alpha\vk}{\cal P}(d\alpha)dx=+\infty. 
\end{equation} 
\end{lem}

{\it Proof}. We shall prove the lemma only for the case that $x_0\in\cSp\cap\cSm$ and $m_+(x_0)-m_-(x_0)>4\pi\beta_+$, 
since the lemma for the other cases are similarly shown. 
In the following, the proof is divided into four steps. \\

{\it Step 1}. 
Fix $x_0\in\cSp\cap\cSm$, and assume that \eqref{eqn:lim-rv} is false to prove the lemma by contradiction. 
Since $\cS\subset\Omega$ now, there exists $0<r_0\ll 1$ such that $B_{2r_0}\Subset\Omega$ and $B_{2r_0}\cap\cS=\{x_0\}$. 
We may assume that $x_0=0$. 
In the following, we consider the problem in $B_{2r_0}$, so that 
\[
\displaystyle 
-\Delta\vk=\lk \int_I \frac{\alpha e^{\alpha\vk}}{\iint_{I\times\Omega} e^{\alpha' \vk} {\cal P}(d\alpha')dx} {\cal P}(d\alpha) \quad\mbox{in $B_{2r_0}$}. 
\]

By retaking $0<r_0\ll 1$, we can take $0<\varepsilon\ll 1$ such that  
\begin{align}
&\frac{(\alpha_+^\ast-2\varepsilon)(\alpha_+^\ast-\varepsilon)}{\varepsilon}\leq \frac{2\pi}{\|s_-\|_{L^1(B_{r_0})}},\quad 
\alpha_+^\ast>2\varepsilon, 
 \label{eqn:lem:judge-rv-1} \\
&(\alpha_+^\ast-2\varepsilon)(m_+(x_0)-m_-(x_0))>4\pi, 
 \label{eqn:lem:judge-rv-2}
\end{align}
since $s_-\in L^1(\Omega)$ and $m_+(0)-m_-(0)>4\pi\beta_+$, recall \eqref{eqn:mq4}. 

Carefully reading \cite{rz2012} shows that there exists $v\in H_{loc}^1(\Omega\setminus\cS)$ such that 
\begin{equation}\label{eqn:lem:judge-rv-3}
\vk\rightarrow v\quad \mbox{in $H_{loc}^1(\Omega\setminus\cS)$}
\end{equation}
and 
\begin{equation}\label{eqn:lem:judge-rv-4}
\left\{
\begin{array}{ll}
\displaystyle 
-\Delta v=(s_+-s_-)+\sum_{y_0\in\cSp}m_+(y_0)\delta_{y_0}-\sum_{y_0\in\cSm}m_-(y_0)\delta_{y_0} & \mbox{in $\Omega$} \\
v=0 & \mbox{on $\partial\Omega$}. 
\end{array}
\right.
\end{equation}

{\it Step 2}. 
Let $z=z(x)$ be the very weak solution of 
\begin{equation*}
\left\{
\begin{array}{ll}
\displaystyle 
-\Delta z=-s_-+(m_+(x_0)-m_-(x_0))\delta_0 & \mbox{in $B_{r_0}$} \\
\displaystyle 
z=b_0:=\min_{\partial B_{r_0}}v & \mbox{on $\partial B_{r_0}$}. 
\end{array}
\right., 
\end{equation*}
see \cite{veron} for the concept of very weak solutions. 
Since $v=v(x)$ satisfies
\[
-\Delta v=(s_+-s_-)+(m_+(x_0)-m_-(x_0))\delta_0\quad \mbox{in $B_{r_0}$},
\]
the maximum principle and $s_+\geq 0$ imply 
\begin{equation}\label{eqn:lem:judge-rv-8}
z\leq v\quad \mbox{a.e. in $B_{r_0}$}.
\end{equation}
Furthermore, we decompose $z=z(x)$ as $z=z_1+z_2$, where $z_1$ and $z_2$ are the solutions of 
\begin{align*}
&-\Delta z_1=(m_+(x_0)-m_-(x_0))\delta_0\ \mbox{in $B_{r_0}$},\quad z_1=b_0\ \mbox{on $\partial B_{r_0}$}, \\
&-\Delta z_2=-s_-\leq 0\ \mbox{in $B_{r_0}$},\quad z_2=0\ \mbox{on $\partial B_{r_0}$}.  
\end{align*}
A direct calculation shows 
\begin{equation}
z_1(x)=\frac{m_+(x_0)-m_-(x_0)}{2\pi}\log\frac{1}{|x|}+c_1,  
 \label{eqn:lem:judge-rv-9}
\end{equation}
for some constant $c_1$ depending only on $b_0$ and $r_0$.\\ 


{\it Step 3}. 
We put 
\[
f=e^{(\alpha_+^\ast-\varepsilon)z_1}\geq 0,\quad 
g=e^{(\alpha_+^\ast-\varepsilon)z_2}\geq 0, \quad  
p=\frac{\alpha_+^\ast-\varepsilon}{\alpha_+^\ast-2\varepsilon}>1. 
\]
Then, \eqref{eqn:lem:judge-rv-9} and \eqref{eqn:lem:judge-rv-2} imply 
\begin{align}
&\left(\int_{B_{r_0/2}}f^{1/p}\right)^p
=\left(\int_{B_{r_0/2}}e^{(\alpha_+^\ast-2\varepsilon)z_1}\right)^p \nonumber\\ 
&\geq a\left(\int_{B_{r_0/2}}|x|^{-\frac{(\alpha_+^\ast-2\varepsilon)(m_+(x_0)-m_-(x_0))}{2\pi}}\right)^p
=+\infty
 \label{eqn:lem:judge-rv-12}
\end{align}
for some $a>0$. 
Noting that $z_2\leq 0$ a.e. in $B_{r_0}$ by the maximum principle, 
we use \eqref{eqn:lem:judge-rv-1} and Lemma \ref{lem:bm-ineq} to obtain 
\begin{align}
&\left(\int_{B_{r_0/2}}g^{-\frac{1}{p-1}}\right)^{-(p-1)}
=\left[\int_{B_{r_0/2}}\exp\left(
\frac{(\alpha_+^\ast-2\varepsilon)(\alpha_+^\ast-\varepsilon)}{\varepsilon}|z_2|
\right)\right]^{-(p-1)} \nonumber\\ 
&\geq \left[\int_{B_{r_0/2}}\exp\left(
\frac{2\pi}{\|s_-\|_{L^1(B_{r_0})}}|z_2|
\right)\right]^{-(p-1)}
\geq (2\pi r_0^2)^{-(p-1)}.
 \label{eqn:lem:judge-rv-13}
\end{align}

{\it Step 4}. We organize \eqref{eqn:lem:judge-rv-12}-\eqref{eqn:lem:judge-rv-13}, the H\"{o}lder inequality, 
\eqref{eqn:lem:judge-rv-8}, \eqref{eqn:lem:judge-rv-3}, the Fatou lemma, $0<r_0\ll 1$ and the assumption of contradiction, 
so that  
\begin{align*}
&+\infty=\left(\int_{B_{r_0/2}}f^{1/p}\right)^p \left(\int_{B_{r_0/2}}g^{-\frac{1}{p-1}}\right)^{-(p-1)} \\
&\leq \int_{B_{r_0/2}}fg
=\int_{B_{r_0/2}}e^{(\alpha_+^\ast-\varepsilon)z}
\leq \int_{B_{r_0/2}} e^{(\alpha_+^\ast-\varepsilon)v}
\leq \liminf_{k\rightarrow\infty}\int_{B_{r_0/2}} e^{(\alpha_+^\ast-\varepsilon)\vk} \\
&\leq \frac{1}{{\cal P}([\alpha_+^\ast-\varepsilon/2,\alpha_+^\ast])}
\liminf_{k\rightarrow\infty}\int_{[\alpha_+^\ast-\varepsilon/2,\alpha_+^\ast]} 
|B_{r_0/2}|^{\frac{\alpha-(\alpha_+^\ast-\varepsilon)}{\alpha}} 
\left(\int_{B_{r_0/2}}e^{\alpha\vk}\right)^{\frac{\alpha_+^\ast-\varepsilon}{\alpha}}{\cal P}(d\alpha) \\
&\leq \frac{1}{{\cal P}([\alpha_+^\ast-\varepsilon/2,\alpha_+^\ast])}
\liminf_{k\rightarrow\infty}\left(1+\iint_{I\times\Omega} e^{\alpha' \vk}\right)<+\infty, 
\end{align*}
a contradiction. \hfill\qed\ \\

We now arrive at the stage to prove \eqref{eqn:residual-vanishing}. \\ 

\underline{{\it Proof of \eqref{eqn:residual-vanishing}}} 
Since $\vk$ is locally uniformly bounded in $\Omega\setminus\cS$, 
and since \eqref{eqn:residual-vanishing} is equivalent to $s_\pm=0$ in \eqref{eqn:rz1}, 
the proof of \eqref{eqn:residual-vanishing} is reduced to showing \eqref{eqn:lim-rv}. 
Moreover, it is reduced to showing \eqref{eqn:pf-rv-1} for $x_0\in\cS$ by virtue of Lemma \ref{lem:judge-rv}.

Property \eqref{eqn:pf-rv-1} is clear for $x_0\in\cSpm\setminus\cSmp$ by \eqref{eqn:mr} and \eqref{eqn:rz5}. 
For $x_0\in\cSp\cap\cSm$, we introduce the sets
\begin{align*}
&{\cal C}=\{(s,t) \ | \ (s-t)^2=8\pi(\beta_+ s+\beta_- t),\ s\geq 0,\ t\geq 0\}, \\
&{\cal D}_+=\{(s,t) \ | \ s-t>4\pi\beta_+,\ s\geq 0,\ t\geq 0 \}, \\
&{\cal D}_{-}=\{(s,t) \ | \ t-s>4\pi\beta_-,\ s\geq 0,\ t\geq 0 \}, 
\end{align*}
see the figure below. 
Then, an elementary calculation shows that ${\cal C}\subset{\cal D}_+\cup{\cal D}_-$, 
which implies \eqref{eqn:pf-rv-1} by \eqref{eqn:mr}. 
The proof is complete. \hfill\qed

\begin{center}
\scalebox{0.75}{\input{neri_fig.tex}}
\end{center}
\section{Proof of Theorem \ref{thm:lbp}}\label{sec:location}

Let $(\lk,\vk)$ be a solution sequence of \eqref{eqn:neri} and assume that the alternative {\bf (II)} in Theorem \ref{thm:main} occurs.

For the purpose, we introcuce 
\[
f_k(t)=\kappa_k\int_I \alpha e^{\alpha t}{\cal P}(d\alpha),\quad 
\kappa_k=\frac{\lk}{\iint_{I\times\Omega}e^{\alpha' v_k}},\quad 
F_k(t)=\kappa_k\int_I e^{\alpha t}{\cal P}(d\alpha). 
\]
Then, \eqref{eqn:neri} reads 
\[
-\Delta v_k=f_k(v_k)\ \mbox{in $\Omega$} \quad v_k=0\ \mbox{on $\partial\Omega$}, 
\]
and we find the following property by combining \eqref{eqn:rz3}-\eqref{eqn:rz4} and \eqref{eqn:residual-vanishing}-\eqref{eqn:mr}:
\begin{equation}\label{eqn:Fk-conv}
F_k(v_k) \overset{\ast}{\rightharpoonup} 
\sum_{x_0'\in\cS}(\beta_+(x_0')m_+(x_0')+\beta_-(x_0')m_-(x_0'))\delta_{x_0'} \quad \mbox{in ${\cal M}(\bar{\Omega})$}, 
\end{equation}
where $\beta_\pm(x_0)$ is as in \eqref{eqn:mq4}. \\

We now comply the complex analysis argument developed by \cite{ye} to prove Theorem \ref{thm:lbp}.\\ 

\underline{{\it Proof of Theorem \ref{thm:lbp}}} 
Fix $x_0\in\cS$ and take $\delta>0$ such that 
\begin{equation}\label{eqn:delta}
B_{2\delta}\subset\Omega,\quad B_{2\delta}\cap\cS=\{x_0\}. 
\end{equation}

Define 
\[
\Gamma=\frac{1}{4\pi}\log(z\bar{z}),\quad 
I_k=\frac{1}{2}(\partial_z\vk)^2,\quad
J_k=\partial_z\Gamma\ast\{\chi_{B_\delta}(\partial_z F_k(\vk))\}
\]
in the usual complex notation $z=X_1+\imath X_2$, 
where $\partial_z=\partial/\partial z$ and $\ast$ denotes the usual convolution. 
Then, we can easily check that 
\[
\partial_{\bar{z}}S_k=0 \ \mbox{in $B_\delta$},\quad 
\mbox{where $S_k=I_k+J_k$ and $\partial_{\bar{z}}=\partial/\partial \bar{z}$}, 
\]
namely, $S_k$ is a holomorphic function. 
Hence, there exists $S_0$, which is holomorphic in $B_\delta$, such that 
\begin{equation}\label{eqn:Sk-conv}
S_k\rightarrow S_0\quad \mbox{locally uniformly in $B_\delta$}. 
\end{equation}

Consider
\begin{equation}\label{eqn:omega}
\omega(x)=(m_+(x_0)-m_-(x_0))H(x,x_0)
+\sum_{x_0'\in\cS\setminus\{x_0\}}(m_+(x_0')-m_-(x_0'))G(x,x_0'). 
\end{equation}
Note that $\omega$ is smooth in $B_\delta$ by \eqref{eqn:delta}. 
Since 
\[
v_k\rightarrow v_0=-(m_+(x_0)-m_-(x_0))\Gamma+\omega \quad \mbox{in $C_{loc}^2(B_\delta\setminus\{0\})$}
\] 
by {\bf (II)-(ii)} in Theorem \ref{thm:main}, it holds that 
\begin{equation}\label{eqn:Ik-conv}
I_k\rightarrow 
I_0=\frac{(m_+(x_0)-m_-(x_0))^2}{32\pi^2z^2}
-\frac{m_+(x_0)-m_-(x_0)}{4\pi z}\omega_z+\frac{\omega_z^2}{2}
\end{equation}
locally uniformly in $B_\delta\setminus\{0\}$. 
Moreover, since $J_k$ takes the another form
\[
J_k=\partial_{zz}\Gamma\ast(\chi_{B_\delta}F_k(\vk))-\partial_z\Gamma\ast\{\partial_z(\chi_{B_\delta})F_k(\vk)\},
\]
\eqref{eqn:Fk-conv} assures that 
\begin{equation}\label{eqn:Jk-conv}
J_k\rightarrow J_0=-\frac{\beta_+m_+(x_0)+\beta_-m_-(x_0)}{4\pi z^2}+J_0'
\end{equation}
locally uniformly in $B_\delta\setminus\{0\}$, where $J_0'$ is the non-singular function defined in $B_\delta$. 

Organizing \eqref{eqn:Sk-conv} and \eqref{eqn:Ik-conv}-\eqref{eqn:Jk-conv}, 
and comparing the coefficients of the singular parts, 
we obtain 
\begin{equation}\label{eqn:coef-inverse-z}
(m_+(x_0)-m_-(x_0))\omega_z(0)=0,\ \mbox{or}\ \omega_z(0)=0
\end{equation}
by $m_+(x_0)-m_-(x_0)\neq 0$ for $x_0\in\cS$. 
Finally, \eqref{eqn:lbp-dom} follows from \eqref{eqn:omega} and \eqref{eqn:coef-inverse-z}. \hfill\qed
\section{Problems on manifolds}\label{sec:mfd}

In this section, we study the Neri mean field equation on manifolds: 
\begin{equation}\label{eqn:neri-mfd}
\left\{
\begin{array}{ll}
\displaystyle 
-\Delta v=\lambda \int_I \frac{\alpha \left(e^{\alpha v}-\frac{1}{|\Omega|}\int_\Omega e^{\alpha v}\right)}{\iint_{I\times \Omega} e^{\alpha' v} {\cal P}(d\alpha')dx} {\cal P}(d\alpha) 
& \mbox{on $\Omega$} \\
\displaystyle 
\int_\Omega vdx=0, &
\end{array}
\right.
\end{equation}
where $v$ is the stream function, 
$\lambda>0$ a constant related to the inverse temperature, 
$\Omega=(\Omega, g)$ a compact and orientable Riemannian surface in dimension two without boundary, 
$g$ the metric on $\Omega$, 
$\Delta=\Delta_g$ the Laplace-Beltrami operator, 
$dx$ the volume element on $\Omega$, 
$|\Omega|$ the volume of $\Omega$, 
${\cal P}\in{\cal M}(I)$ a Borel probability measure on $I$, 
${\cal M}(I)$ the space of measures on $I$, 
and $I=[-1,1]$. 

To state the results, we prepare some notations. 
Let $(\lk,\vk)$ be a solution sequence of \eqref{eqn:neri-mfd}. 
Similarly to Section \ref{sec:intro}, we define the blowup set $\cS$ by 
\begin{align*}
&\cS=\cSp\cup\cSm, \\
&\cS_{\pm}=\{x_0\in\Omega\ | \ \mbox{there exists $\xk\in\Omega$ such that  $\xk\rightarrow x_0$ and $\vk(\xk)\rightarrow\pm\infty$} \}, 
\end{align*}
and introduce the measure functions $\nu_{k,\pm}=\nu_{k,\pm}(dx)\in{\cal M}(\Omega)$ and $\mu_k=\mu_k(d\alpha dx)\in{\cal M}(I\times\Omega)$ defined by 
\begin{align}
&\nu_{k,\pm}
=\lk\int_{I_\pm}\frac{|\alpha| e^{\alpha \vk}}{\iint_{I\times\Omega}e^{\alpha' \vk}}{\cal P}(d\alpha), 
 \label{eqn:def-nuk-mfd}\\
&\mu_k(d\alpha dx)=\lk\frac{e^{\alpha \vk}}{\iint_{I\times\Omega}e^{\alpha' \vk}}{\cal P}(d\alpha)dx, 
\end{align}
respectively, where $I_+=[0,1]$ and $I_-=[-1,0]$. 

With this notation, we review the result of \cite{rz2012}.

\begin{pro}\label{fact:rz-mfd}
Let $(\lk,\vk)$ be a solution sequence of \eqref{eqn:neri-mfd} with $\lk\rightarrow\lzero$ for some $\lzero\geq 0$. 
Then, passing to a subsequence, we have the following alternatives. 

{\bf (I')} Compactness: $\displaystyle \limsup_{k\rightarrow\infty}\|\vk\|_\infty<+\infty$, that is, $\cS=\emptyset$.

\noindent
Then, there exists $v\in E$ such that 
$\vk\rightarrow v$ in $E$ and $v$ is a solution of \eqref{eqn:neri} for $\lambda=\lzero$, where 
\[
E=\{v\in H^1(\Omega)\ | \ \int_\Omega v=0\}. 
\]

{\bf (II')} Concentration:  $\displaystyle \limsup_{k\rightarrow\infty}\|\vk\|_\infty=+\infty$, that is, $\cS\neq\emptyset$.

\noindent
Then, $\cS=\cSp\cup\cSm$ is finite and the following properties {\bf a')}-{\bf c')} hold: 

{\bf a')} There exists $0\leq s_\pm\in L^1(\Omega)\cap L_{loc}^\infty(\Omega\setminus\cS)$ such that
\begin{equation}\label{eqn:rz1-mfd}
\nu_{k,\pm} \overset{\ast}{\rightharpoonup} \nu_\pm=s_\pm+\sum_{x_0\in\cS_\pm}m_\pm(x_0)\delta_{x_0} \quad \mbox{in ${\cal M}(\Omega)$}, 
\end{equation}
with $m_\pm(x_0)\geq 4\pi$ for every $x_0\in\cS_\pm$, 
where $\delta_{x}\in{\cal M}(\Omega)$, $x\in\Omega$, denotes the Dirac measure centered at $x$. 

{\bf b')} There exist $\zeta_{x_0}\in{\cal M}(I)$ and $0\leq r\in L^1(I\times\Omega)$ such that 
\begin{equation}\label{eqn:rz2-mfd}
\mu_k \overset{\ast}{\rightharpoonup} 
\mu=\mu(d\alpha dx)
=r(\alpha,x){\cal P}(d\alpha)dx+\sum_{x_0\in\cS}\zeta_{x_0}(d\alpha)\delta_{x_0}(dx) \quad \mbox{in ${\cal M}(I\times\Omega)$}. 
\end{equation}

{\bf c')} For every $x_0\in\cS$, we have  
\begin{align}
&8\pi\int_I \zeta_{x_0}(d\alpha)=\left(\int_I \alpha \zeta_{x_0}(d\alpha)\right)^2
 \label{eqn:rz3-mfd}\\
&m_\pm(x_0)=\int_{I_\pm}|\alpha| \zeta_{x_0}(d\alpha),\quad s_\pm(x)=\int_{I_\pm}|\alpha| r(\alpha,x){\cal P}(d\alpha), 
 \label{eqn:rz4-mfd}
\end{align}
where $m_\pm(x_0)$ is as in \eqref{eqn:rz1-mfd}. 
Moreover, for every $x_0 \in\cS_\pm\setminus\cS_\mp$, it holds that
\begin{equation}\label{eqn:rz5-mfd}
m_\mp(x_0)=\int_{I_\mp}|\alpha|\zeta_{x_0}(d\alpha)=0. 
\end{equation}
\end{pro}

As already stated in Section \ref{sec:intro}, 
we shall show the results corresponding to Theorems \ref{thm:main} and \ref{thm:lbp}, 
except for {\bf (II)-(i)} in Theorem \ref{thm:main}.
The first result is

\begin{thm}\label{thm:main-mfd}
Assume that the alternative {\bf (II')} in Proposition \ref{fact:rz-mfd} occurs. 
Then, it holds that $v_k \rightarrow v_0$ in $C_{loc}^2(\Omega\setminus\cS)$, where
\begin{equation}\label{eqn:residual-vanishing-mfd}
v_0(x)=\sum_{x_0'\in\cS}(m_+(x_0')-m_-(x_0'))G(x,x_0')
\end{equation}
with $m_\pm(x_0')\geq 4\pi$ for every $x_0'\in\cS$, 
and where $G=G(x,y)$ is the Green function defined by 
\begin{equation}\label{eqn:green-mfd}
-\Delta G(x,y)=\delta_y-\frac{1}{|\Omega|}\ \mbox{on $\Omega$},\quad \int_\Omega G(x,y)dx=0. 
\end{equation}
Moreover, we have
\begin{equation}\label{eqn:mq3-mfd}
(m_+(x_0)-m_-(x_0))^2=8\pi(\beta_+(x_0)m_+(x_0)+\beta_-(x_0)m_-(x_0))
\end{equation}
for every $x_0\in\cS$, where $\beta_\pm(x_0)$ is as in \eqref{eqn:mq4}. 
\end{thm}

To state the second result, we introduce several notations. 

Let $(\lk,\vk)$ be a solution sequence of \eqref{eqn:neri-mfd}. 
Given $x_0\in\Omega$, we take an iso-thermal chart $(\Psi,U)$ such that 
\begin{equation}\label{eqn:iso-thermal}
\Phi(x_0)=0,\quad g=e^{\xi(X)}(dX_1^2+dX_2^2),\quad X=\Psi(x),\quad \xi(0)=0. 
\end{equation}
Then, $\tvk(X)=\vk\circ\Psi^{-1}(X)$ is the solution of 
\begin{equation}\label{eqn:tvk}
-\Delta_X \tvk=e^\xi \left(
\lk \int_I \frac{\alpha \left(e^{\alpha \tvk}-\frac{1}{|\Omega|}\int_\Omega e^{\alpha \vk}\right)}{\iint_{I\times\Omega} e^{\alpha' \vk}} {\cal P}(d\alpha) 
\right)\quad \mbox{in $\Psi(U)$}.  
\end{equation}
Furthermore, we introduce the regular part $H_\Psi=H_\Psi(x,y)$ of the Green function $G=G(x,y)$, 
defined by \eqref{eqn:green-mfd}, 
relative to the isothermal chart satisfying \eqref{eqn:iso-thermal}, namely, 
\begin{equation}\label{eqn:robin-mfd}
H_\Psi(x,y)=G(x,y)-\frac{1}{2\pi}\log\frac{1}{|X-Y|}
\end{equation}
for $(x,y)\in U\times U$, where $X=\Psi(x)$ and $Y=\Psi(y)$. 

Under these preparations, the location of the blowup points is characterized in terms of $G$, $H_\Psi$ and $\xi$ as follows. 

\begin{thm}\label{thm:lbp-mfd}
Let $(\lk,\vk)$ be a solution sequence of \eqref{eqn:neri-mfd} 
and assume that the alternative {\bf (II')} in Proposition \ref{fact:rz-mfd} occurs. 
Then, for every $x_0\in{\cal S}$, we have 
\begin{align}
&\nabla_X\left[
H_\Psi(\Psi^{-1}(X),x_0)+\frac{\xi(X)}{8\pi}
\right.\nonumber\\ 
&\qquad \left.\left.
+\sum_{x_0'\in{\cal S}\setminus\{x_0\}}\frac{m_+(x_0')-m_-(x_0')}{m_+(x_0)-m_-(x_0)}G(\Psi^{-1}(X),x_0')
\right]\right|_{X=0}=0. 
 \label{eqn:lbp-mfd}
\end{align}
\end{thm}\ \\

The remainder of this section is devoted to showing Theorems \ref{thm:main-mfd} and \ref{thm:lbp-mfd}. 
The proofs are analogue to those of Theorems \ref{thm:main} and \ref{thm:lbp}, 
and therefore we shall give only a sketch of them here. \\

\underline{{\it Proof of Theorem \ref{thm:main-mfd}}} 
At first, \eqref{eqn:mq3-mfd} is shown in the same way as \eqref{eqn:mr}, see Section \ref{sec:rv}. 

Next, we shall give a sketch of the proof of the following property 
corresponding to Lemma \ref{lem:judge-rv} in Section \ref{sec:rv}: 
for $x_0\in\cS$, if 
\[
m_+(x_0)-m_-(x_0)>4\pi\beta_+\quad \mbox{or}\quad m_-(x_0)-m_+(x_0)>4\pi\beta_-
\] 
then 
\begin{equation}\label{eqn:judge-mfd}
\lim_{k\rightarrow\infty}\iint_{I\times\Omega}e^{\alpha\vk}{\cal P}(d\alpha)dx=+\infty. 
\end{equation}

Here, we shall only investigate the case that $x_0\in\cSp\cap\cSm$ and $m_+(x_0)-m_-(x_0)>4\pi\beta_+$, 
since property \eqref{eqn:judge-mfd} is similarly shown for the other cases. 

Fix $x_0\in\cSp\cap\cSm$, assume that \eqref{eqn:judge-mfd} fails, and take an isothermal chart $(U,\psi)$ around $x_0$ as above. 
Then, $\vk=\vk(\xi^{-1}(X))$ satisfies 
\[
-\Delta_X \vk=e^{\xi}\lambda_k
\int_{I}\frac{\alpha e^{\alpha\vk}}{\iint_{I\times\Omega}e^{\alpha' \vk}}{\cal P}(d\alpha)\quad 
\mbox{in $\psi(U)$}
\]
and there exists $0<r_0\ll 1$ such that $B_{2r_0}\subset\psi(U)$ and $\psi^{-1}(B_{2r_0})\cap\cS=\{x_0\}$. 
For simplicity, we shall write $X$ and $\Delta_X$ by $x$ and $\Delta$, respectively. 
In the following, we consider the problem on $B_{2r_0}$, so that 
\[
-\Delta \vk=e^{\xi}\lambda_k
\int_{I}\frac{\alpha e^{\alpha \vk}}{\iint_{I\times\Omega}e^{\alpha' \vk}}{\cal P}(d\alpha)\quad 
\mbox{in $B_{2r_0}$}. 
\]

By retaking $0<r_0\ll 1$, we obtain $0<\varepsilon\ll 1$ such that  
\begin{align*}
&\frac{(\alpha_+^\ast-2\varepsilon)(\alpha_+^\ast-\varepsilon)}{\varepsilon}\leq \frac{2\pi}{\|s_-\|_{L^1(B_{r_0})}},\quad 
\alpha_+^\ast>2\varepsilon, \\
&(\alpha_+^\ast-2\varepsilon)(m_+(x_0)-m_-(x_0))>4\pi. 
\end{align*}

Here, we note that there exist $v\in H_{loc}^1(\Omega\setminus\cS)$ and $c_0\in{\bf R}$ such that 
\[
\vk\rightarrow v\quad \mbox{in $H_{loc}^1(\Omega\setminus\cS)$}
\]
and 
\[
\displaystyle 
-\Delta v=(s_+-s_-)+\sum_{y_0\in\cSp}m_+(y_0)\delta_{y_0}-\sum_{y_0\in\cSm}m_-(y_0)\delta_{y_0}-c_0\quad \mbox{on $\Omega$}
\]
with $\int_\Omega v=0$, see \cite{rz2012}.

Let $z=z(x)$ be the very weak solution of
\[
\left\{
\begin{array}{ll}
\displaystyle 
-\Delta z=-s_-+(m_+(x_0)-m_-(x_0))\delta_0-c_0 & \mbox{in $B_{r_0}$} \\
\displaystyle 
z=b_0:=\min_{\partial B_{r_0}}v & \mbox{on $\partial B_{r_0}$}, 
\end{array}
\right.
\]
see \cite{veron} for the concept of very weak solutions. 
Note that $v=v(x)$ satisfies 
\[
-\Delta v=(s_+-s_-)+(m_+(x_0)-m_-(x_0))\delta_0-c_0\quad \mbox{in $B_{r_0}$}.
\]
Moreover, we decompose $z=z(x)$ as $z=z_1+z_2$, where $z_1$ and $z_2$ are the solutions of 
\begin{align*}
&-\Delta z_1=(m_+(x_0)-m_-(x_0))\delta_0-c_0\ \mbox{in $B_{r_0}$},\quad z_1=b_0\ \mbox{on $\partial B_{r_0}$}, \\
&-\Delta z_2=-s_-\leq 0\ \mbox{in $B_{r_0}$},\quad z_2=0\ \mbox{on $\partial B_{r_0}$},
\end{align*}
respectively. 
Putting
\[
f=e^{(\alpha_+^\ast-\varepsilon)z_1}\geq 0,\quad 
g=e^{(\alpha_+^\ast-\varepsilon)z_2}\geq 0, \quad  
p=\frac{\alpha_+^\ast-\varepsilon}{\alpha_+^\ast-2\varepsilon}>1, 
\]
one can show the following estimates similarly to Section \ref{sec:rv}: 
\[
\left(\int_{B_{r_0/2}}f^{1/p}\right)^p=+\infty,\quad 
\left(\int_{B_{r_0/2}}g^{-\frac{1}{p-1}}\right)^{-(p-1)}\geq (2\pi r_0^2)^{-(p-1)}.
\]
Similarly to Section \ref{sec:rv} again, we obtain 
\[
+\infty=\left(\int_{B_{r_0/2}}f^{1/p}\right)^p \left(\int_{B_{r_0/2}}g^{-\frac{1}{p-1}}\right)^{-(p-1)}<+\infty, 
\]
a contradiction. 
Hence, \eqref{eqn:judge-mfd} is established. 

Finally, \eqref{eqn:residual-vanishing-mfd} is shown by the argument developed in the last part of Section \ref{sec:rv}, 
and the proof is complete. \hfill\qed\\ 

\underline{{\it Proof of Theorem \ref{thm:lbp-mfd}}} 
Let $(\lk,\vk)$ be a solution sequence of \eqref{eqn:neri-mfd} 
and assume that the alternative {\bf (II')} in Proposition \ref{fact:rz-mfd} occurs. 

We put 
\[
f_k(t)=\kappa_k\int_I \alpha e^{\alpha t}{\cal P}(d\alpha)+c_k,\quad 
\kappa_k=\frac{\lk}{\iint_{I\times\Omega}e^{\alpha' v_k}},\quad 
c_k=-\frac{\lk \iint_{I\times\Omega}\alpha e^{\alpha v_k}}{|\Omega|\iint_{I\times\Omega}e^{\alpha v_k}}. 
\]
Then \eqref{eqn:neri-mfd} is equivalent to 
\[
-\Delta v_k=f_k(v_k)\ \mbox{on $\Omega$} \quad \mbox{with}\ \int_\Omega v_k=0. 
\]
Furthermore, similarly to Section \ref{sec:location}, we obtain 
\begin{equation}\label{eqn:Fk-conv'}
F_k(v_k) \overset{\ast}{\rightharpoonup} 
\sum_{x_0'\in\cS}(\beta_+(x_0')m_+(x_0')+\beta_-(x_0')m_-(x_0'))\delta_{x_0'}
+c_0v_0 \quad \mbox{in ${\cal M}(\Omega)$}, 
\end{equation}
where $\beta_\pm(x_0')$ is as in \eqref{eqn:mq4}, 
$v_0$ as in \eqref{eqn:residual-vanishing-mfd}, $c_0=\lim_{k\rightarrow\infty}c_k$ and 
\[
F_k(t)=\kappa_k\int_I e^{\alpha t}{\cal P}(d\alpha)+c_k t. 
\]

To end the proof of Theorem \ref{thm:lbp-mfd}, we again comply the complex analysis argument developed by \cite{ye}.\\ 

Fix $x_0\in\cS$ and take an iso-thermal chart satisfying \eqref{eqn:iso-thermal}. 
We write $\tvk(X)=\vk\circ\Psi^{-1}(X)$ by $\vk(X)$ for simplicity. 
Then \eqref{eqn:tvk} reads 
\begin{equation}\label{eqn:vk-C'}
-\Delta_X \vk=e^\xi f_k(\vk)\quad \mbox{in $\Psi(U)$}, 
\end{equation}
and there exists $\delta>0$ such that 
\begin{equation}\label{eqn:delta'}
B_{2\delta}\subset \Psi(U),\quad \Psi^{-1}(B_{2\delta})\cap\cS=\{x_0\}. 
\end{equation}

Putting 
\[
\Gamma=\frac{1}{4\pi}\log(z\bar{z}),\quad 
I_k=\frac{1}{2}(\partial_z\vk)^2,\quad
J_k=\partial_z\Gamma\ast\{e^\xi \chi_{B_\delta}(\partial_z F_k(\vk))\}, 
\]
we calculate 
\[
\partial_{\bar{z}}S_k=0 \ \mbox{in $B_\delta$},\quad 
\mbox{where $S_k=I_k+J_k$ and $\partial_{\bar{z}}=\partial/\partial \bar{z}$}, 
\]
and hence there exists $S_0$, which is holomorphic in $B_\delta$, such that 
\begin{equation}\label{eqn:Sk-conv'}
S_k\rightarrow S_0\quad \mbox{locally uniformly in $B_\delta$}. 
\end{equation}
Furthermore, we consider 
\begin{equation}\label{eqn:omega'}
\omega(X)=(m_+(x_0)-m_-(x_0))H_\Psi(x,x_0)
+\sum_{x_0'\in\cS\setminus\{x_0\}}(m_+(x_0')-m_-(x_0'))G(x,x_0'), 
\end{equation}
recall $X=\Psi(x)$, 
where $H_\Psi=H_\Psi(x,y)$ is the regular part of the Green function defined by \eqref{eqn:robin-mfd}. 
Note that $\omega$ is smooth in $B_\delta$ by the smoothness of $\Omega$ and \eqref{eqn:delta'}. 
Since $v_k\rightarrow v_0=-(m_+(x_0)-m_-(x_0))\Gamma+\omega$ in $C_{loc}^2(B_\delta\setminus\{0\})$ 
by \eqref{eqn:residual-vanishing-mfd}, it holds that 
\begin{equation}\label{eqn:Ik-conv'}
I_k\rightarrow 
I_0=\frac{(m_+(x_0)-m_-(x_0))^2}{32\pi^2z^2}
-\frac{m_+(x_0)-m_-(x_0)}{4\pi z}\omega_z+\frac{\omega_z^2}{2}
\end{equation}
locally uniformly in $B_\delta\setminus\{0\}$. 
Also, since $J_k$ takes the another form
\[
J_k=\partial_{zz}\Gamma\ast(e^\xi\chi_{B_\delta}F_k(\vk))-\partial_z\Gamma\ast\{\partial_z(e^\xi\chi_{B_\delta})F_k(\vk)\},
\]
\eqref{eqn:Fk-conv'} and $\xi(0)=0$ imply that 
\begin{align}
J_k\rightarrow J_0=&-\frac{\beta_+(x_0)m_+(x_0)+\beta_-(x_0)m_-(x_0)}{4\pi z^2} \nonumber\\
&-\frac{\beta_+(x_0)m_+(x_0)+\beta_-(x_0)m_-(x_0)}{4\pi z}\xi_z(0)+J_0'
 \label{eqn:Jk-conv'}
\end{align}
locally uniformly in $B_\delta\setminus\{0\}$, where $J_0'$ is the non-singular function defined in $B_\delta$. 

Organizing \eqref{eqn:Sk-conv'} and \eqref{eqn:Ik-conv'}-\eqref{eqn:Jk-conv'}, 
and comparing the coefficients of the singular parts, 
we obtain 
\begin{equation}\label{eqn:coef-inverse-z'}
\frac{m_+(x_0)-m_-(x_0)}{\beta_+(x_0)m_+(x_0)+\beta_-(x_0)m_-(x_0)}\omega_z(0)+\xi_z(0)=0. 
\end{equation}
Consequently, \eqref{eqn:lbp-mfd} follows from \eqref{eqn:coef-inverse-z'}, \eqref{eqn:omega'} and \eqref{eqn:mq3-mfd}.
\hfill\qed

\section*{Acknowledgements}
This research is partially supported 
by {\it Programma di scambi internazionali con universit\`{a} ed istituti di ricerca stranieri per la mobilit\`{a} di breve durata di docenti, 
ricercatori e studiosi} of Universit\`{a} di Napoli Federico II 
and by {\it Engineering Science Young Researcher Dispatch Program} of Osaka University.

\end{document}